\documentclass[11pt]{article}

\usepackage{amsmath}
\usepackage{amsfonts}
\usepackage{color}
\usepackage{hyperref}
\usepackage{multirow}
\usepackage{graphicx}
\usepackage{enumitem}

\newtheorem{teo}{Theorem}[section]
\newtheorem{coro}{Corollary}[section]
\newtheorem{defi}{Definition}[section]
\newtheorem{assump}{Assumption}

\newcommand{\halmos}{\hfill$\Box$}
\newenvironment{pro}{\noindent\textit{Proof:}}{\halmos}

\newcommand{\R}{\mathbb{R}}
\newcommand{\xtrial}{x^{\mathrm{trial}}} 
\newcommand{\half}{\frac{1}{2}}
\newcommand{\ftarget}{f_{\mathrm{target}}}
\newcommand{\Mbar}{\underline{M}_{\bar{x}}}
\newcommand{\f}{\underline{f}}
\newcommand{\n}{\underline{n}}
\newcommand{\g}{\underline{g}}
\newcommand{\T}{\underline{T}}
\newcommand{\Pe}{\underline{P}}
\newcommand{\Ome}{\underline{\Omega}}
\newcommand{\ele}{\underline{\ell}}
\newcommand{\uu}{\underline{u}}
\newcommand{\h}[1]{\textcolor{blue}{\underline{#1}}}
\newcommand{\mbar}{\bar{m}}
\DeclareMathOperator*{\Minimize}{Minimize}
\DeclareMathOperator*{\argmin}{argmin}

\begin{document}

\title{On complexity and convergence of high-order\\ coordinate
  descent algorithms for smooth nonconvex box-constrained
  minimization\footnote{This work was supported by FAPESP (grants
    2013/07375-0, 2016/01860-1, and 2018/24293-0) and CNPq (grants
    302538/2019-4 and 302682/2019-8).}}

\author{
  V. S. Amaral\thanks{Dept. of Applied Mathematics,
    Institute of Mathematics, Statistics, and Scientific Computing,
    University of Campinas, 13083-859, Campinas, SP,
    Brazil. email: vitalianoamaral@hotmail.com, andreani@ime.unicamp.br,
    martinez@ime.unicamp.br}
  \and
  R. Andreani\footnotemark[2]
  \and
  E. G. Birgin\thanks{Dept. of Computer Science, Institute of
    Mathematics and Statistics, University of S\~ao Paulo, Rua do
    Mat\~ao, 1010, Cidade Universit\'aria, 05508-090, S\~ao Paulo, SP,
    Brazil. email: egbirgin@ime.usp.br}
  \and
  D. S. Marcondes\thanks{Dept. of Applied Mathematics, Institute of
    Mathematics and Statistics, University of S\~ao Paulo, Rua do
    Mat\~ao, 1010, Cidade Universit\'aria, 05508-090, S\~ao Paulo, SP,
    Brazil. email: diaulas@ime.usp.br}
  \and
  J. M. Mart\'{\i}nez\footnotemark[2]
}

\date{December 20, 2021\footnote{This version is a revision of the
    version submitted in August 3, 2021.}}

\maketitle

\begin{abstract}
Coordinate descent methods have considerable impact in global
optimization because global (or, at least, almost global) minimization
is affordable for low-dimensional problems. Coordinate descent methods
with high-order regularized models for smooth nonconvex
box-constrained minimization are introduced in this work. High-order
stationarity asymptotic convergence and first-order stationarity
worst-case evaluation complexity bounds are established. The computer
work that is necessary for obtaining first-order
$\varepsilon$-stationarity with respect to the variables of each
coordinate-descent block is $O(\varepsilon^{-(p+1)/p})$ whereas the
computer work for getting first-order $\varepsilon$-stationarity with
respect to all the variables simultaneously is
$O(\varepsilon^{-(p+1)})$.  Numerical examples involving
multidimensional scaling problems are presented. The numerical
performance of the methods is enhanced by means of coordinate-descent
strategies for choosing initial points.\\
  

\end{abstract}

\section{Introduction} \label{introduction}

In order to minimize a multivariate function it is natural to keep
fixed some of the variables and to modify the remaining ones trying to
decrease the objective function value. Coordinate descent (CD) methods
proceed systematically in this way and, many times, obtain nice
approximations to minimizers of practical optimization problems.
Wright~\cite{wright} surveyed traditional approaches and modern
advances on the introduction and analysis of CD methods. Although the
CD idea is perhaps the most natural one to optimize functions, it
received little attention from researchers due to poor performance in
many cases and lack of challenges in terms of convergence
theory~\cite{powell}. The situation changed dramatically in the last
decades. CD methods proved to be useful for solving machine learning,
deep learning and statistical learning problems in which the number of
variables is big and the accuracy required at the solution is moderate
\cite{boyd,nesterov}. Many applications arose and, in present days,
efficient implementations and insightful theory for understanding the
CD properties are the subject of intense research. See, for example,
\cite{abrs,bt,bst,bone,bs,bh,cd,ey,llx,xu,ywbm} among many others.

In this paper we are concerned with complexity issues of CD methods
that employ high-order models to approximate the subproblems that
arise at each iteration. The use of high-order models for
unconstrained optimization was defined and analyzed from the point of
view of worst-case complexity in \cite{bgmst1} and subsequent papers
\cite{bgms,cagtholder,grap1,grap2,jmholder,novo3}. In \cite{bgms}
numerical implementations with quartic regularization were introduced.
In \cite{cagtholder}, \cite{grap1}, \cite{grap2}, and~\cite{jmholder},
new high-order regularization methods were introduced with H\"older,
instead of Lipschitz, conditions on the highest-order derivatives
employed.  In \cite{novo2}, high-order methods were studied as
discretizations of ordinary differential equations. These methods
generalize the methods based on third-order models introduced in \cite
{griewank} and later developed in
\cite{cagtpart1,cagtpart2,crs,dussault,nesterovpolyak} among many
others. Griewank \cite{griewank} introduced third-order regularization
having in mind affine scaling properties. Nesterov and Polyak
\cite{nesterovpolyak} introduced the first cubic regularized Newton
methods with better complexity results than the ones that were known
for gradient-like algorithms \cite{gyy}. In \cite{novo1}, a multilevel
strategy that exploits a hierarchy of problems of decreasing dimension
was introduced in order to reduce the global cost of the step
computation. However, high-order methods remain difficult to implement
in the many-variables case due to the necessity of computing
high-order derivatives and solving nontrivial model-based
subproblems. Nevertheless, if the number of variables is small,
high-order model-based methods are reliable alternatives to classical
methods. This feature can be exploited in the CD framework.

High-order models are interesting from the point of view of global
optimization because, many times, local algorithms get stuck at points
that satisfy low-order optimality conditions from which one is able to
escape using high-order resources. The escaping procedure is
affordable if one restricts the search to low-dimensional subspaces,
which suggests the employment of CD procedures.
   
This paper is organized as follows. In Section~\ref{background}, we
present some background on optimality conditions, while in
Section~\ref{preliminar} we survey a high-order algorithmic framework
that provides a basis for the development of CD algorithms. In
Section~\ref{cda}, we present block CD methods that, for each
approximate minimization on a group of variables, employ high-order
regularized subproblems and we prove asymptotic convergence. In
Section~\ref{complexity} we prove worst-case complexity results. In
Section~\ref{discussion} the obtained theoretical results are
discussed. In Section~\ref{exp}, we study a family of problems for
which CD is suitable and we include a CD-strategy that improves
convergence to global solutions. Conclusions are given in
Section~\ref{concl}.\\

\noindent
\textbf{Notation.} The symbol $\|\cdot\|$ denotes the Euclidean norm.

\section{Background on high-order optimality conditions} \label{background}

In order to understand the main results of this paper we need to visit
the topic of necessary optimality conditions of high order. The main
question is: Which is the relation between minimizers of a function
and minimizers of its Taylor polynomials?  Firstly, we show that, in
one variable, the two concepts are closely related in the sense that
local minimizers of a function are local minimizers of all its Taylor
polynomials. Immediately, we show with a simple counterexample that
this property is not true if the number of variables is greater
than~1. The third step is to show that, for an arbitrary number of
variables, every minimizer of $f$ is a minimizer of its Taylor
polynomials regularized by a suitable Lipschitz constant. This
definition leads us to distinguish between exclusive and inclusive
optimality conditions. Exclusive conditions are the ones that can be
expressed exclusively in terms of the function derivatives. Inclusive
ones are related with a slightly more global behavior and include
Lipschitz bounds. Inclusive conditions are stronger than exclusive
ones. In this paper, we show that algorithmic limit points are more
related to inclusive conditions than to exclusive ones.

As it is well known from elementary calculus, if a function $\f:\R \to
\R$ possesses derivatives up to order~$p$ at~$\bar x \in \R$,
denoted by $\f^{(j)}$ for $j=1,\dots,p$, its Taylor
polynomial of order~$p$ around~$\bar x$ is given by
\[
\T_p(\bar x, x) = \f(\bar x) + \sum_{j=1}^p \frac{1}{j!} \f^{(j)}(\bar
x) (x - \bar x)^j.
\]
If~$\f$ and its derivatives up to order~$p$ are continuous and
$\f^{(p)}$ satisfies a Lipschitz condition defined by $\gamma_1 > 0$
in a neighborhood of $\bar x$, we know that
\begin{equation} \label{ftp}
|\f(x) - \T_p(\bar x, x)| \leq \frac{\gamma_1}{(p+1)!} |x - \bar x|^{p+1}
\end{equation}
for all~$x$ in a neighborhood of~$\bar x$. This fact allows one to
prove the necessary optimality condition given in
Theorems~\ref{background1} and~\ref{background2}.

\begin{teo}  \label{background1}
Assume that $\f: \R \to \R$, its derivatives up to order $p$ are
continuous, and $\f^{(p)}$ satisfies a Lipschitz condition defined by
$\gamma_1 > 0$ in a neighborhood of $x^*$. Assume, moreover, that $a <
b$, $x^*$ is a local minimizer of $\f$ subject to $x \in [a,b]$, and
there exists $q \leq p$ such that $\f^{(j)}(x^*) = 0$ for
$j=1,\dots,q-1$ and $\f^{(q)}(x^*) \neq 0$. Then,
\begin{enumerate}
\item if $q$ is even, then we have that $\f^{(q)}(x^*) > 0$;
\item if $a < x < b$, then $q$ is even;
\item if $x = a$ and $q$ is odd, then $\f^{(q)}(x^*) > 0$;
\item if $x = b$ and $q$ is odd, then $\f^{(q)}(x^*) < 0$.
\end{enumerate}
\end{teo}

\begin{pro}
Suppose that $q \leq p$ is such that all the derivatives of order $j <
q \leq p$ are null and $\f^{(q)}(x^*) \neq 0$. Then, by (\ref{ftp}),
\[
\left|\f(x) -\f(x^*) - \left[ \frac{1}{q!}\f^{(q)}(x^*)(x - x^*)^q
  + \dots + \frac{1}{p!} \f^{(p)}(x^*)(x - x^*)^p \right] \right|
\leq \frac{\gamma_1}{(p+1)!} |x - x^*|^{p+1}.
\]
Then,
\[
\resizebox{\textwidth}{!}{
$\left|\f(x) - \f(x^*) - \frac{1}{q!}\f^{(q)}(x^*)(x - x^*)^q
\right| - \left| \frac{1}{(q+1)!}\f^{(q+1)}(x^*)(x - x^*)^{q+1}
\dots + \frac{1}{p!} \f^{(p)}(x^*)(x - x^*)^p \right| \leq
\frac{\gamma_1}{(p+1)!} |x - x^*|^{p+1}$.}
\]
Thus, if $p=q$, it follows trivially that
\begin{equation} \label{paradivi}
\left| \f(x) - \f(x^*) - \frac{1}{q!}\f^{(q)}(x^*)(x - x^*)^q
\right| \leq c |x - x^*|^{q+1}.
\end{equation}  
If $p > q$, for all $j=q+1,\ldots,p$, the quantities
$|\frac{1}{j!}f^{(j)}(x^*)|$ are bounded by the same constant. By the
boundedness of $|x-x^*|$ in a neighborhood of $x^*$ and the fact that
$p+1 > q+1$, (\ref{paradivi}) follows as well.
%Thus, with an appropriate restriction of the neighborhood that defines
%the local minimization of $\f$, and using continuity of the
%derivatives, there exists $c > 0$ such that
%\begin{equation} \label{paradivi}
%\left| \f(x) - \f(x^*) - \frac{1}{q!}\f^{(q)}(x^*)(x - x^*)^q
%\right| \leq c |x - x^*|^{q+1}.
%\end{equation}  
Assume firstly that $q$ is even. Then, dividing (\ref{paradivi}) by
$(x - x^*)^q > 0$, we have that
\begin{equation} \label{para2}
\left| \frac{ \f(x)-\f(x^*)}{(x - x^*)^q} -
\frac{1}{q!}\f^{(q)}(x^*) \right| \leq c |x - x^*|.
\end{equation}          
Taking limits for $x \to  x^*$ we deduce that
\begin{equation} \label{para3}
\lim_{x \to x^*} \left| \frac{ \f(x)-\f(x^*)}{(x - x^*)^q} -
\frac{1}{q!}\f^{(q)}(x^*) \right| = 0.
\end{equation}
Thus, 
\begin{equation} \label{therhs}
\lim_{x \to x^*} \frac{ \f(x)-\f(x^*)}{(x - x^*)^q} =
\frac{1}{q!} \f^{(q)}(x^*).
\end{equation}  
Since $\f(x) \geq \f(x^*)$ for all $x$ sufficiently close to $x^*$ and
the right-hand side of (\ref{therhs}) is different from zero, we
deduce that $ \f^{(q)}(x^*) > 0$.  Therefore, we proved that if not
all the derivatives are null, the first statement in the thesis is
true.

Now consider the case in which all the derivatives of order $j < q
\leq p$ are null, $a < x^* < b$, and $f^{(q)}(x^*) \neq 0$. Suppose,
by contradiction that $q$ is odd. Assume, firstly, that
$x>x^*$. Dividing (\ref{paradivi}) by $(x - x^*)^q > 0$, we have that
\begin{equation} \label{repe1}
\left| \frac{ \f(x)-\f(x^*)}{(x - x^*)^q} -
\frac{1}{q!} \f^{(q)}(x^*) \right| \leq c |x - x^*|.
\end{equation}   
Taking lateral limits for $x > x^*$ and  $x \to  x^*$ we deduce that
\begin{equation} \label{repe2}
\lim_{x \to x^*, \; x > x^*} \left| \frac{ \f(x)-\f(x^*)}{(x -
  x^*)^q} - \frac{1}{q!} \f^{(q)}(x^*) \right| = 0.
\end{equation}
Thus, 
\begin{equation} \label{therhs2}
\lim_{x \to x^*, \; x > x^*} \frac{ \f(x)-\f(x^*)}{(x - \bar
  x)^q} = \frac{1}{q!} \f^{(q)}(x^*).
\end{equation}  
Since $f(x) \geq f(x^*)$ for all $x$ sufficiently close to $ x^*$, we
deduce that $ \f^{(q)}(x^*) \geq 0$. A similar reasoning for $x < x^*$
leads to $ \f^{(q)}(x^*) \leq 0$. Therefore, $ \f^{(q)}(x^*) = 0$.
Therefore, we proved that if all the derivatives of order $j < q \leq
p$ are null, $a < x < b$, and $ \f^{(q)}(x^*) \neq 0$, then $q$ is
even.

Let us prove now that, if all the derivatives of order $j < q \leq p$
are null, $\f^{(q)}(x^*) \neq 0$, $x^* = a$ and $q$ is odd, we have
that $\f^{(q)}(x^*) > 0$. Dividing (\ref{paradivi}) by $(x - x^*)^q >
0$, we obtain (\ref{repe1}), (\ref{repe2}), and (\ref{therhs2}) with
$x^* = a$. Since $ \f(x) \geq \f(x^*)$ for all $x$ sufficiently close
to $x^*$ and, by assumption, $ \f^{(q)}(x^*) \neq 0$, we have that $
\f^{(q)}(x^*) \geq 0$.  The last part of the thesis follows exactly in
the same way.
\end{pro}

\begin{teo} \label{background2}
Assume that $\f: \R \to \R$ and its derivatives up to order $p$ are
continuous and~$\f^{(p)}$ satisfies a Lipschitz condition defined by
$\gamma_1 > 0$ in a neighborhood of~$x^*$. Assume, moreover,
that~$x^*$ is a local minimizer of~$\f$. Then, $x^*$ is a local
minimizer of the Taylor polynomial $\T_p(x^*, x)$.
\end{teo}

\begin{pro}
By Theorem~\ref{background1} we have four alternatives for the
coefficients of the Taylor polynomial of order $p$. The first one is
that all its coefficients are null. In this case, $x^*$ is, trivially,
a minimizer of the polynomial and there is nothing to prove.

In the second case the first nonnull coefficient of the polynomial is
positive and its order is even. Therefore, the Taylor polynomial can
be written as
\[
\T_p(x^*, x) = \f(x^*) + \sum_{j=q}^p \frac{1}{j!} \f^{(j)}(
x^*) (x - x^*)^j
\]   
for some even $q\leq p$ and $\frac{1}{q!} \f^{(j)}(x^*)>0$. Then,
\begin{equation} \label{follow}
\frac{\T_p(x^*, x) - \f(x^*)}{(x-x^*)^q} = \frac{1}{q!}
\f^{(q)}(x^*) + \sum_{j=q+1}^p \frac{1}{j!} \f^{(j)}(x^*) (x -
x^*)^{j-q}.
\end{equation}    
This implies that $x^*$ is a local minimizer of $\T_p(x^*, x)$ as we
wanted to prove.

In the third case $x^* = a$, $q$ is odd and $ \frac{1}{q!}
\f^{(j)}(x^*)>0$.  Then, (\ref{follow}) takes place and $a$ is a local
minimizer. The fourth case, in which $x^* = b$ and $ \frac{1}{q!}
\f^{(j)}(x^*)<0$, follows in a similar way.
\end{pro}\\

We now consider the $\n$-dimensional case. If $\f:\R^{\n} \to \R$
admits continuous derivatives up to order $p \in \{1, 2, 3, \dots\}$,
then the Taylor polynomial of order $p$ of $\f$ around $x^*$ is
defined as
\begin{equation} \label{taylorp}
\T_p(x^*, x) = \f(x^*) + \sum_{j=1}^p \Pe_j(x^*, x),
\end{equation}
where ${\Pe}_j(x^*, x)$ is an homogeneous polynomial of degree $j$
given by
\begin{equation} \label{Pj}
{\Pe}_j(x^*, x) = \frac{1}{j!} \left( (x_1 - x^*_1)
\frac{\partial}{\partial x_1} + \dots + (x_n - x^*_n)
\frac{\partial}{\partial x_n} \right)^j \f(x).
\end{equation}
For completeness we define ${\Pe}_0(x^*, x) = \f(x^*)$.

Let us define $\varphi(t) = \f (x^* + t (x - x^*))$. Obviously, if
$x^*$ is a local minimizer of $\f$ over a nonempty closed and convex
set $C \subset \R^{\n}$, it turns out that $0$ is a local minimizer of
$\varphi(t)$ for every choice of $x \in C$. Thus, by
Theorem~\ref{background2}, $0$ is a local minimizer of the Taylor
polynomial associated with $\varphi$ subject to the interval defined
by the boundary of $C$.  But, by the construction of (\ref{taylorp}),
this implies that $x^*$ is a minimizer of $\T_p(x^*, x)$ along any
line that passes through $x^*$ over the interval defined by the
boundary of $C$.  This fact is stated in Theorem~\ref{background3}.

\begin{teo} \label{background3}
Assume that $\f: \R^{\n} \to \R$ and its derivatives up to order $p$
are continuous and satisfy a Lipschitz condition in a neighborhood of
$x^*$. Assume, moreover, that $x^*$ is a local minimizer of $\f$. Let
$\mathcal{L}$ be a line that passes through $x^*$.  Then, $x^*$ is a
local minimizer of the Taylor polynomial $\T_p(x^*, x)$ subject to
$\mathcal{L} \cap C$.
\end{teo}

\begin{pro}
Observe that the fact that the derivatives of order $p$ satisfy a
Lipschitz condition imply that the $p$-th derivative of $\varphi$
exhibits the same property. Then, apply Theorem~\ref{background2}.
\end{pro}

\begin{defi} \label{def31}
We say that $x^*$ is $p$th-order stationary of $\f$ over the closed
and convex set $C$ if, for all $x \in C$, $0$~is a local minimizer of
the Taylor polynomial of order $p$ that corresponds to the univariate
function $\varphi(t) = \f(x^* + t(x - x^*))$ restricted to the
constraint $x^* + t (x - x^*) \in C$.
\end{defi}        

\noindent
\textbf{Counterexample.} Unfortunately, it is not true that, when
$x^*$ is a local minimizer of $\f$, it is also a local minimizer of
the associated Taylor polynomial. (As we saw in
Theorem~\ref{background2}, this property is indeed true when $\n =
1$.) For example, if $\f(x_1, x_2) = x_2^2 - x_1^2 x_2 + x_1^4$, we
have that $(0,0)$ is a global minimizer of $\f$, but it is not a local
minimizer of its Taylor polynomial of order $p=3$.\\

%Suppose that $x^* \in C$ is a local minimizer of~$\f$ over the closed
%and closed and convex set~$C$ and all the $p$-th derivatives of~$\f$
%satisfy a Lipschitz condition in a ball centered in~$x^*$. By standard
%multidimensional Calculus, we know that there exists $\gamma > 0$ such
%that for all $x \in C$ in that ball,
%\begin{equation} \label{tay2}
%\f(x) \leq \T_p(x^*, x) + \gamma \|x - x^*\|^{p+1}.
%\end{equation}
%The following theorem defines a necessary optimality condition for the
%minimization of~$f$ onto~$C$.

%\begin{teo} \label{background3bis}
%Assume that $\f: \R^{\n} \to \R$ and its derivatives up to order $p$
%are continuous and satisfy a Lipschitz condition in a neighborhood of
%$x^*$ in such a way that (\ref{tay2}) holds. Assume, moreover, that
%$x^*$ is a local minimizer of $\f$ onto the closed and convex set $C$.
%Assume that $\sigma \geq \gamma$. Then, $x^*$ is a local minimizer of
%$\T_p(x^*, x) + \sigma \|x - x^*\|^{p+1}$ subject to $x \in C$.
%\end{teo}

%\begin{pro}
%If the derivatives of $\f$ up to order $p$ are continuous, local
%minimizers are necessarily $p$th-order stationary according to
%Definition~\ref{def31}.
%\end{pro}

%The following theorem motivates a stronger definition of
%$p$-stationarity that holds when a Lipschitz condition is satisfied.

In the following theorem we prove that, although according to the
counterexample above, a minimizer does not need to minimize the Taylor
polynomial, such property is true if the Taylor polynomial is
regularized with a Lipschitz term.

\begin{teo} \label{arbitrary}
Assume that ${\cal D} \subset \R^{\n}$, $\f:{\cal D} \to \R$,
and~$x^*$ is a local minimizer of~$\f(x)$ over~${\cal D}$ such that,
for all $x \in {\cal D}$,
\begin{equation} \label{lili}
\f(x) \leq \T_p(x^*, x) + \gamma  \|x - x^*\|^{p+1},
\end{equation}
where $\T_p$ is, as defined in~(\ref{taylorp}), the Taylor polynomial
of order $p$ of $\f$ around~$x^*$. Then, for all $\sigma \geq \gamma$,
$x^*$ is a local minimizer of $\T_p(x^*, x) + \sigma \|x -
x^*\|^{p+1}$ over ${\cal D}$.
\end{teo}
   
\begin{pro} 
Suppose that the thesis is not true. Then, $x^*$ is not a local
minimizer of $\T_p(x^*, x) + \gamma \|x - x^*\|^{p+1}$ over ${\cal
  D}$. Thus, there exists $\{x^k\} \subset {\cal D}$ such that
$\lim_{k \to \infty} x^k = x^*$ and
\[
\T_p(x^*, x^k) + \gamma \|x^k - x^*\|^{p+1} < \T_p(x^*,x^*) = \f(x^*).
\]
Thus, by (\ref{lili}), 
\[
\f(x^k) < \f(x^*)
\]  
for all $k = 0, 1, 2, \dots$ This contradicts the fact that $x^*$ is
a local minimizer of $\f$ over ${\cal D}$.
\end{pro}\\

The following definition is motivated by Theorem~\ref{arbitrary}.

\begin{defi} \label{def32}
Assume that ${\cal D} \subset \R^{\n}$, $\f:{\cal D} \to \R$, $x^*$ is
such that (\ref{lili}) holds for all $x \in {\cal D}$, and that
$\sigma \geq \gamma$. Then $x^* \in {\cal D}$ is said to be
$p$th-order $\sigma$-stationary of $\f$ over ${\cal D}$ if $x^*$ is a
local minimizer of $\T_p(x^*, x) + \sigma \|x - x^*\|^{p+1}$ over
${\cal D}$.
\end{defi}

\begin{figure}[ht!]
\begin{center}
\input{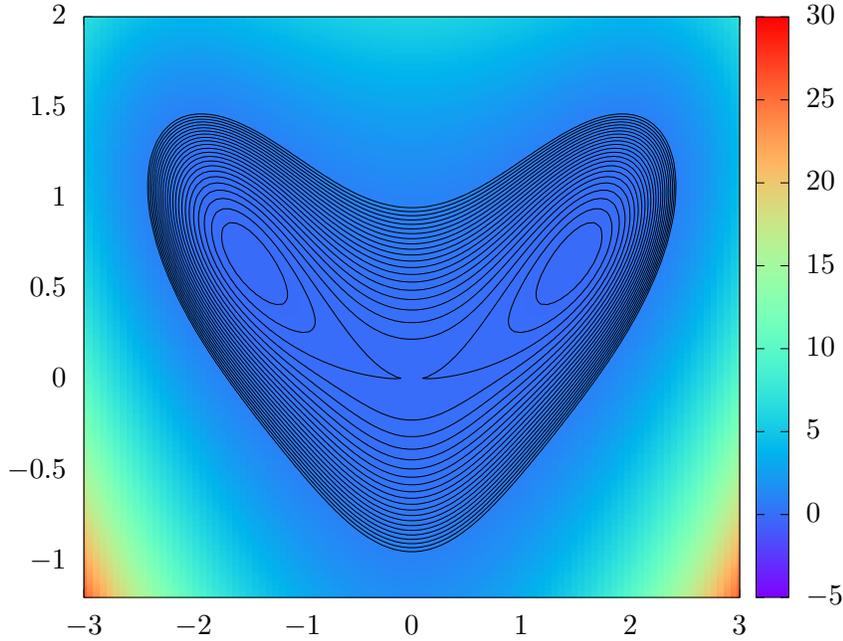}
\end{center}
\caption{Level sets of $\T_p((0,0),(x_1,x_2)) + \sigma \| (x_1,x_2) -
  (0,0) \|^{p+1}$ with $p=3$ and $\sigma=0.125$, where
  $\T_p((0,0),(x_1,x_2))$ is the $p$th-order Taylor polynomial of
  $\f(x_1, x_2) = x_2^2 - x_1^2 x_2 $ (that coincides with $\f$). The
  graphic shows that Condition {\bf C5}
  with $p=3$ and $\sigma = 0.125$ does not hold at $(0,0)$, since it
  is \textit{not} a local minimizer of the regularized $p$th-order
  Taylor polynomial. There are two local minimizers at ``the eyes of
  the cat''.}
\label{fig1}
\end{figure}

It is trivial to see that, if ${\cal D}$ is convex and $x^*$ is
$p$th-order $\sigma$-stationary of $\f$ over ${\cal D}$ according to
Definition~\ref{def32}, then it is $p$th-order $\tilde
\sigma$-stationary for every $\tilde \sigma \geq \sigma$ and it is
also $p$th-order stationary according to
Definition~\ref{def31}. However, $p$th-order $\sigma$-stationarity is
strictly stronger than $p$th-order stationarity. Consider the function
$\f(x_1, x_2) = x_2^2 - x_1^2 x_2 $ and $p = 3$. Note that $x^* = (0,
0)$ satisfies (\ref{lili}) with $\gamma = 0$. Straightforward
calculations show that the point $(0,0)$, that is not a local
minimizer of~$\f$, is $p$th-order stationary according to
Definition~\ref{def31}. On the other hand, $(0,0)$ is not $p$th-order
$\sigma$-stationarity if $\sigma < 1/4$. See Figure~\ref{fig1}.

At this point it is convenient to summarize the properties of
candidates to solutions of $\mbox{Minimize } \f(x) \mbox{ subject to }
x \in C$, where $C$ is closed and convex. Let us consider the
following conditions with respect to $x^* \in C$:
\begin{description}
\item[\textbf{C1:}] $x^*$ is a local minimizer.
\item[\textbf{C2:}] $x^*$ is a local minimizer of the Taylor polynomial
  over every feasible segment that passes through $x^*$.
\item[\textbf{C3:}] $x^*$ is a local minimizer of the Taylor polynomial
  around $x^*$.
\item[\textbf{C4:}] $x^*$ is a local minimizer of $\T_p(x^*, x) + \gamma
  \|x - x^*\|^{p+1}$, where $\gamma$ is a Lipschitz constant.
\item[\textbf{C5:}] $x^*$ is a local minimizer of $\T_p(x^*, x) + \sigma
  \|x - x^*\|^{p+1}$, where $\sigma > \gamma$ and $\gamma $ is a
  Lipschitz constant.
\item[\textbf{C6:}] $x^*$ is a local minimizer of $\T_p(x^*, x) + \sigma
  \|x - x^*\|^{p+1}$, where $0 < \sigma < \gamma$ and $\gamma $ is a
  Lipschitz constant.
\end{description}  
We proved that \textbf{C1}, \textbf{C2}, \textbf{C4} and \textbf{C5}
are necessary optimality conditions, while \textbf{C3} and \textbf{C6}
are not. We also showed that \textbf{C1} $\Rightarrow$ \textbf{C4}
$\Rightarrow$ \textbf{C5}, and \textbf{C3} $\Rightarrow$ \textbf{C6}
$\Rightarrow$ \textbf{C4} $\Rightarrow$ \textbf{C5}. However,
\textbf{C1} does not imply neither \textbf{C3} nor \textbf{C6}.

\begin{defi} \label{def33}
We say that an optimality condition is {\it exclusive} if it can be
verified using only values of the derivatives up to order $p$ at the
point under consideration.
\end{defi}

Optimality conditions that are not exclusive are said to be
\textit{inclusive}. Only condition \textbf{C2} above is
exclusive. \textbf{C4} and \textbf{C5} are inclusive necessary
optimality conditions because they use information on the Lipschitz
constant in a neighborhood of $x^*$. Thus, the information that they
require is not restricted to derivatives of order at most $p$ at a
single point. The annihilation of the gradient at $x^*$ and the
positive semidefiniteness of the Hessian are exclusive first-order and
second-order necessary optimality conditions for unconstrained
optimization. The most natural high-order exclusive optimality
condition for convex constrained optimization is
\textbf{C2}. In~\cite{cagt5}, an exclusive optimality condition based
on curves was presented. However, exclusive necessary optimality
conditions are essentially weaker than inclusive ones. In fact, assume
that $x^*$ satisfies \textbf{C5} and that \textbf{C} is an arbitrary
exclusive necessary optimality condition. Then, $x^*$ is a local
minimizer of $ \T_p(x^*, x) + \sigma \|x - x^*\|^{p+1}$, where $\sigma
> \gamma$ and $\gamma $ is a Lipschitz constant. Then, $x^*$ satisfies
the exclusive condition {\bf C} for the minimization of $\T_p(x^*, x)
+ \sigma \|x - x^*\|^{p+1}$. Then, since \textbf{C} is a necessary
optimality condition, it is satisfied by $x^*$ for the local
minimization of $\T_p(x^*, x) + \sigma \|x - x^*\|^{p+1}$.  But all
the derivatives up to order $p$ of $\T_p(x^*, x) + \sigma \|x -
x^*\|^{p+1}$ exist at $x^*$ and coincide with the derivatives up to
order $p$ of $\f$. So, $x^*$ satisfies \textbf{C} for the minimization
of $\f$.

In order to see that \textbf{C5} is strictly stronger than \textbf{C}
(for every exclusive necessary optimality condition \textbf{C}),
consider the functions $\f(x_1, x_2) = x_2^2 - x_1^2 x_2 $ and $F(x_1,
x_2) = x_2^2 - x_1^2 x_2 + x_1^4$. The origin $x^*=(0, 0)$ is a local
(and global) minimizer of $F$, therefore, it must satisfy the
necessary exclusive optimality condition \textbf{C} of order
$p=3$. Since, up to order $p=3$, the derivatives of $\f$ and $F$ are
the same, it turns out that $x^*$ satisfies the necessary optimality
condition \textbf{C} of order $p=3$, applied to the minimization of
$\f$.  (Note that $x^*$ is not a local minimizer of $\f$.) However,
$x^*$ does not satisfy condition \textbf{C5} if $\sigma < 1/4$. In
this case, every $\sigma > 0$ is bigger than the Lipschitz constant of
$\f$ associated with third-order derivatives, thus, we found an
example in which the exclusive condition \textbf{C} holds but the
inclusive condition \textbf{C5} does not.
   
\section{Regularized high-order minimization with box constraints} \label{preliminar}

In this section, we consider the problem
\begin{equation} \label{mainprob}
\mbox{Minimize } \f(x) \mbox{ subject to } x \in \Ome,
\end{equation}
where $\Ome \subset \R^{\n}$ is given by
\begin{equation} \label{omegaunderline}
\Ome = \{x \in \R^{\n} \;|\; \ele \leq x \leq \uu \}
\end{equation}
and $\ele, \uu \in \R^{\n}$ are such that $\ele < \uu$. We assume that
$\f$ has continuous first derivatives into $\Ome$. We denote $\g(x) =
\nabla \f(x)$ and $\g_P(x) = P_{\Ome}(x - \g(x)) - x$, for all $x \in
\Ome$, where $P_{\Ome}$ is the Euclidean projection operator
onto~$\Ome$. In the remaining of this section, the results
from~\cite{bmnewgen} that are relevant to the present work are
surveyed and a natural extension of the main algorithm
in~\cite{bmnewgen}, that makes it possible to consider a wider class
of models, is introduced.

Each iteration~$k$ of Algorithm~2.1 introduced in \cite{bmnewgen}
computes a new iterate~$x^{k+1}$ satisfying $(p+1)$th-order descent
with respect to $\f(x^k)$ through the approximate minimization of a
$(p+1)$th-regularized model of the function~$\f$ around
the iterate~$x^k$. For all $\bar x \in \R^n$, let $\Mbar: \R^n \to \R$
be a ``model'' of $\f(x)$ around $\bar x$; and assume that $\nabla
\Mbar(x)$ exists for all $x \in \Ome$. We now present an algorithm
that corresponds to a single iteration of the algorithm introduced
in~\cite{bmnewgen}.\\

\noindent
\textbf{Algorithm~\ref{preliminar}.1.} Assume that $p \in \{ 1, 2, 3,
\dots\}$, $\alpha > 0$, $\sigma_{\min} > 0$, $\tau_2 \geq \tau_1 >1$,
$\theta > 0$, and $\bar x \in \Ome$ are given.

\begin{description}
\item[Step 1.] Set $\sigma \leftarrow 0$.

\item[Step 2.] Compute $\xtrial \in \Ome$ such that
  \begin{equation} \label{modelbaja}
  \Mbar(\xtrial) + \sigma \|\xtrial - \bar x\|^{p+1} \leq   \Mbar(\bar x) 
  \end{equation}
  and 
  \begin{equation} \label{paradasub}
  \left\| P_{\Ome}\left[\xtrial - \left. \nabla \left( \Mbar(x) + 
  \sigma \|x - \bar x\|^{p+1} \right) \right|_{x=\xtrial} \right] - \xtrial  \right\|
  \leq \theta \|\xtrial - \bar x\|^p.
  \end{equation}              

\item[Step 3.] If
  \begin{equation} \label{armijo4}
  \f(\xtrial) \leq \f(\bar x) - \alpha \| \xtrial - \bar x \|^{p+1},
  \end{equation}
  then define $x^{+} = \xtrial$ and stop returning $x^{+}$ and
  $\sigma$. Otherwise, update $\sigma \leftarrow \max\{ \sigma_{\min},
  \tau \sigma \}$ with $\tau \in [\tau_1, \tau_2]$ and go to Step 2.
\end{description}

\noindent
\textbf{Remark.} The trial point $\xtrial$ computed at Step~2 is
intended to be an approximate solution to the subproblem
\begin{equation} \label{subproblem}
\mbox{Minimize } \Mbar(x) + \sigma \| x - \bar x \|^{p+1}
\mbox{ subject to } x \in \Ome.
\end{equation}
Note that conditions (\ref{modelbaja}) and (\ref{paradasub}) can
always be achieved. In fact, by the compactness of $\Ome$,
if~$\xtrial$ is a global minimizer of~(\ref{subproblem}), then it
satisfies the condition
\[
\left\| P_{\Ome} [\xtrial - \nabla(\Mbar(x) + 
\sigma \|x - \bar x\|^{p+1})\big|_{x=\xtrial}] - \xtrial  \right\| = 0;
\]
and so (\ref{paradasub}) takes place. In addition, if $\xtrial$ is a
global minimizer, since $\bar x$ is a feasible point,
(\ref{modelbaja}) must hold as well.

\begin{assump} \label{a4}
There exists $L > 0$ such that, for all $\xtrial$ computed by
Algorithm~\ref{preliminar}.1, $x = \xtrial$ satisfies
\begin{equation} \label{lips2}
\left\| \g (x) - \nabla \Mbar(x) \right\| \leq L \|x - \bar x\|^p, 
\end{equation}
\begin{equation} \label{lips3}
\Mbar(\bar x) = \f(\bar x) \mbox{ and } \f(x) \leq \Mbar(x) + L \|x -
\bar x\|^{p+1}.
\end{equation} 
\end{assump}

If $\Mbar(x)$ is the Taylor polynomial of order~$p$ of~$\f$
around~$\bar x$ and the $p$th-order derivatives of~$\f$ satisfy a
Lipschitz condition with Lipschitz constant~$L$, then
Assumption~\ref{a4} is satisfied. However, the situations in which
Assumption~\ref{a4} holds are not restricted to the case in which
$\Mbar(x) = \T_p(\bar x, x)$. For example, we may choose $\Mbar(x) =
\f(x)$. (Note that, in this case, $p$ may be arbitrarily large but
only first derivatives of $\f(x)$ need to exist.) Although the results
in \cite{bmnewgen} only mention the choice $\Mbar(x) = \T_p(\bar x,
x)$, these results only depend on Assumption~\ref{a4}. Thus, they can
be trivially extended to the general choice of $\Mbar(x)$.

\begin{teo} \label{teo1}
Suppose that Assumption~\ref{a4} holds. If the regularization
parameter $\sigma$ in~(\ref{modelbaja}) satisfies $\sigma \geq L +
\alpha$, then the trial point $\xtrial$ satisfies the sufficient
descent condition~(\ref{armijo4}). Moreover,
\begin{equation}\label{crucial}
\left\| \g_P (x^{+}) \right\| \leq \left( L + \tau_2 \left( L + \alpha
\right) (p+1) + \theta \right) \|x^{+}-\bar x\|^p
\end{equation}  
and 
\begin{equation} \label{cuantoba}
\f(x^{+}) \leq \f(\bar x) - \alpha \left( \frac{\| \g_P (x^{+}) \|}{L
  + \tau_2 \left(L + \alpha \right) (p+1) + \theta} \right)^{(p+1)/p}.
\end{equation}               
\end{teo}

\begin{pro}
This theorem condensates the results in~\cite[Lemmas~3.2--3.4]{bmnewgen}.
\end{pro}\\

Theorem~\ref{teo1} justifies the definition of an algorithm for
solving (\ref{mainprob}) based on repetitive application of
Algorithm~\ref{preliminar}.1 and shows that such algorithm enjoys good
properties in terms of convergence and complexity. On the one hand,
each iteration of the algorithm requires $O(1)$ functional evaluations
and finishes satisfying a suitable sufficient descent condition. On
the other hand, that condition implies that infinitely many iterations
with gradient-norm bounded away from zero are not possible if the
function is bounded below. Moreover, (\ref{cuantoba}) leads to a
complexity bound on the number of iterations based on the norm of the
projected gradient. In the following sections, we prove that, thanks
to Theorem~\ref{teo1}, similar convergence and evaluation complexity
properties hold for a coordinate descent algorithm.

\section{High-order coordinate descent algorithm} \label{cda}

In this section, we consider the problem
\begin{equation} \label{problema}
\mbox{Minimize } f(x) \mbox{ subject to } x \in \Omega,
\end{equation}
where $\Omega \subset \R^{n}$ is given by
\begin{equation} \label{omega}
\Omega = \{x \in \R^{n} \;|\; \ell \leq x \leq u   \}
\end{equation}
and $\ell, u \in \R^n$ are such that $\ell < u$. We assume that $f$
has continuous first derivatives over $\Omega$.

At each iteration of the coordinate descent method introduced in this
section for solving~(\ref{problema}), (i) a nonempty set of indices
$I_k \subseteq \{1, \dots, n\}$ is selected, (ii) coordinates
corresponding to indices that are not in $I_k$ remain fixed, and (iii)
Algorithm~\ref{preliminar}.1 is applied to the minimization of $f$
over~$\Omega$ with respect to the free variables, i.e.\ variables with
indices in $I_k$. From now on, given $v \in \R^n$, we denote by $v_I
\in \R^{|I|}$ the vector whose components are the components of $v$
whose indices belong to $I \subseteq \{1,\dots,n\}$. For all $x \in
\Omega$, we define $g_{P,I}(x) \in \R^n$ by
\[
[g_{P,I}(x)]_i =
\left\{
\begin{array}{cl}
[g_P(x)]_i, &  \mbox{if } i \in I, \\
0, & \mbox{if } i \notin I.
\end{array}
\right.
\]   
Since $\Omega$ is a box, this definition is equivalent to $g_{P,I}(x)
= P_\Omega(x - g_I(x)) - x$, where
\[
[g_I(x)]_i =
\left\{
\begin{array}{cl}
[g(x)]_i, &  \mbox{if } i \in I, \\
0, & \mbox{if } i \notin I.
\end{array}
\right.
\]   
This equivalence, that will be used in the theoretical convergence
results below, is not true if $\Omega$ is an arbitrary closed and
convex set. This is the reason for which we consider CD algorithms
only with box constraints.\\

\noindent
\textbf{Algorithm~\ref{cda}.1.} Assume that $p \in \{ 1, 2, 3,
\dots\}$, $\alpha > 0$, $\sigma_{\min} > 0$, $\tau_2 \geq \tau_1 >1$,
$\theta > 0$, and $x^0 \in \Omega$ are given. Initialize $k \leftarrow
0$.

\begin{description}
  
\item[Step 1.] Choose a nonempty set $I_k \subseteq \{1, \dots, n\}$.

\item[Step 2.] Consider the problem
  \begin{equation} \label{reduzido}
  \mbox{Minimize } f(x) \mbox{ subject to }
  x \in \Omega \mbox{ and } x_i = x_i^k \mbox{ for all } i \notin I_k.
  \end{equation}
  Let $\bar x = x^k_{I_k}$. Setting $\f$, $\Ome$, and $\Mbar$
  properly, apply Algorithm~\ref{preliminar}.1 to obtain~$x^+$
  and $\sigma_k$.

\item[Step 3.] Define $x^{k+1}$ as $x^{k+1}_{I_k} = x^+$ and
  $x_i^{k+1}=x_i^k$ for all $i \not\in I_k$, set $k \leftarrow k+1$,
  and go to Step~1.

\end{description}

\begin{assump} \label{a5}
There exists $L > 0$ such that for all $k$, $\bar x$, $\f$, and
$\Mbar$ set at the $k$th iteration of Algorithm~\ref{cda}.1 and for
all $\xtrial$ computed by Algorithm~\ref{preliminar}.1 when called at
the $k$th iteration of Algorithm~\ref{cda}.1, (\ref{lips2}) and
(\ref{lips3}) take place with $x=\xtrial$.
\end{assump}

If $\Mbar(x)$ is the Taylor polynomial of order~$p$ of~$\f$
around~$\bar x$ and the $p$th-order derivatives of~$f$ satisfy a
Lipschitz condition with Lipschitz constant~$L$, then
Assumption~\ref{a5} is satisfied.

\begin{teo} \label{teo2}
Suppose that Assumption~\ref{a5} holds. Then, there exists $c > 0$,
which only depends on $L$, $\tau_2$, $\alpha$, $p$, and $\theta$ such
that, for all $k = 0, 1, 2, \dots$, the point $x^{k+1}$ computed by
Algorithm~\ref{cda}.1 is well defined and satisfies
\begin{equation} \label{armijo5}
f(x^{k+1}) \leq f(x^k) - \alpha \| x^{k+1} - x^k \|^{p+1} 
\end{equation}
and
\begin{equation}\label{crucial2}
\left\| g_{P,I_k} (x^{k+1}) \right\| \leq c \|x^{k+1}-x^k\|^p.  
\end{equation}  
\end{teo}
 
\begin{pro}
(\ref{armijo5}) follows from (\ref{armijo4}), while (\ref{crucial2})
  follows from the application of Theorem~\ref{teo1}.
\end{pro}

\begin{teo} \label{teo3}
Suppose that Assumption~\ref{a5} holds. Let $\{x^k\}$ be the sequence
generated by Algorithm~\ref{cda}.1. Then,
\begin{equation} \label{limxk}
\lim_{k \to \infty} \|x^{k+1} - x^k \| = 0, 
\end{equation}
\begin{equation} \label{limgplus}
\lim_{k \to \infty}  \left\| g_{P,I_k} (x^{k+1}) \right\| = 0,
\end{equation}
and
\begin{equation} \label{limgp}
\lim_{k \to \infty}  \left\| g_{P,I_k} (x^k) \right\| = 0.
\end{equation} 
\end{teo}

\begin{pro}
Since~$\Omega$ is compact, we have that~$f$ is bounded below
onto~$\Omega$. Thus, (\ref{limxk}) follows from~(\ref{armijo5}) and,
in consequence, (\ref{limgplus}) follows from~(\ref{limxk}) and
(\ref{crucial2}). Let us prove (\ref{limgp}). Assume that $I \subseteq
\{1, \dots, n\}$ is nonempty and arbitrary.  By the continuity of the
gradient, the function $\|g_{P,I}(x)\|$ is continuous for all $x \in
\Omega$ and, since $\Omega$ is compact, it is uniformly
continuous. Then, given $\varepsilon > 0$, there exists $\delta_I>0$
such that, whenever $\|x - y\| \leq \delta_I$, we have that
$\|g_{P,I}(x) - g_{P,I}(y)\| \leq \varepsilon/2$.  Since the number of
different subsets of $\{1, \dots, n\}$ is finite, we have that $\delta
\equiv \min \{\delta_I\;|\; \emptyset \neq I \subseteq \{1, \dots,
n\}\} > 0$. Thus, for all $I \subseteq \{1, \dots, n\}$, if $\|x - y\|
\leq \delta$, we have that $\|g_{P,I}(x) - g_{P,I}(y)\| \leq
\varepsilon/2$. Now, by (\ref{limxk}), there exists $k_0$ such that,
whenever $k \geq k_0$, we have that $ \|x^{k+1} - x^k \| \leq
\delta$. Then, by the definition of $\delta$, if $k \geq k_0$,
$\|g_{P,I}(x^{k+1}) - g_{P,I}(x^k)\| \leq \varepsilon/2$ for all
nonempty $I \subseteq \{1, \dots, n\}$. In particular, taking $I =
I_k$, if $k \geq k_0$, we have that $\|g_{P,I_k}(x^{k+1}) -
g_{P,I_k}(x^k)\| \leq \varepsilon/2$.  Finally, by (\ref{limgplus}),
there exists $k_1 \geq k_0$ such that, for all $k \geq k_1$,
$\|g_{P,I_k}(x^{k+1})\| \leq \varepsilon/2$. By the triangular
inequality, adding the last two inequalities we have that
$\|g_{P,I_k}(x^{k})\| \leq \varepsilon$. Since $\varepsilon>0$ was
arbitrary, this completes the proof of~(\ref{limgp}).
\end{pro}\\

The following assumption guarantees that all the indices $i \in \{1,
\dots, n\}$ belong to some $I_k$ at least every $\bar m$ iterations.
This guarantees that the CD method tries to reduce the function with
respect to each variable $x_i$ infinitely many times.

\begin{assump} \label{aemebar}
There exists $\bar m < +\infty$ such that, for all $i \in
\{1,\dots,n\}$:
\begin{enumerate}
\item There exists $k \leq \bar m$ such that $i \in I_k$; 
\item For any $k \in \mathbb{N}$, if $i \in I_k$, then there exists $m
  \leq \bar m$ such that $i \in I_{k+m}$.
\end{enumerate}
\end{assump} 

Note that Assumption~\ref{aemebar} allows us to use not only cyclic
versions, but also random versions of the CD method. In particular,
the block of coordinates chosen at each iteration can be chosen at
random, with the condition that, every $\bar{m}$ iterations, all
blocks are chosen at least once.
  
\begin{teo} \label{teo4}
Suppose Assumptions~\ref{a5} and~\ref{aemebar} hold. Let $\{x^k\}$ be the
sequence generated by Algorithm~\ref{cda}.1. Then,
\begin{equation} \label{limgpxk}
\lim_{k \to \infty} \|g_P(x^k) \| = 0. 
\end{equation}
Moreover, if $x^* \in \Omega$ is a limit point of $\{x^k\}$, then we
have that $\|g_P(x^*)\|=0$.
\end{teo}

\begin{pro}
Let $i \in \{1, \dots, n\}$. By Assumption~\ref{aemebar}, there exists
an infinite set of increasing indices $K = \{ k_1, k_2, k_3, \dots\}$
such that $i \in I_{k_\ell}$ and $k_{\ell+1} \leq k_\ell + \bar m$ for
all $\ell=1, 2, 3, \dots$ Then, by (\ref{limgp}) in
Theorem~\ref{teo3}, since, by definition, given $I \subseteq
\{1,\dots,n\}$, $[g_{P,I}(x)]_i = [g_P(x)]_i$ for any $i \in I$,
\begin{equation} \label{gpixk}
\lim_{k \in K} [g_P(x^k)]_i = 0.
\end{equation}
Let $j \in \{1, 2, \dots\}$ be arbitrary. By (\ref{limxk}), the
triangular inequality, and the uniform continuity of $g_P$, we have
that
\[
\lim_{k \in K} |[g_P(x^{k+j})]_i - [g_P(x^k)]_i| = 0.
\]
Therefore, by (\ref{gpixk}), 
\begin{equation} \label{inparti}
\lim_{k \in K}  [g_P(x^{k+j})]_i = 0.
\end{equation}   
In particular, (\ref{inparti}) holds for all $j = 1, \dots, \bar
m$. This implies that
 \begin{equation} \label{inparti2}
\lim_{k \to \infty}  [g_P(x^{k})]_i = 0.
\end{equation}     
 Thus, the thesis is proved.
\end{pro}\\

Theorem~\ref{teo4} shows that limit points of sequences generated by
Algorithm~\ref{cda}.1 are first-order stationary. The rest of this
section is dedicated to prove that, under suitable conditions,
$p$th-order stationarity with respect to each variable also
holds. More precisely, if the same nonempty set $I_k$ is repeated
infinitely many times, $p$-stationarity holds in the limit for the
variables $x_i$ with $i \in I_k$. For this purpose, we need to define
different notions of stationarity.

In Theorem~\ref{teo4} we proved that Algorithm~\ref{cda}.1 is
satisfactory from the point of view of first-order stationarity. In
the CD approach we cannot advocate for full stationarity of high order
because cross derivatives that involve variables that are never
optimized together are not computed at all. However, if optimization
with respect to the same group of variables occurs at infinitely many
iterations, it is reasonable to conjecture that high-order optimality
with respect to those variables would, in the limit, take place. For
obtaining such result, it is not enough to satisfy criteria
(\ref{modelbaja}) and (\ref{paradasub}) when solving subproblems. The
reason is that condition~(\ref{paradasub}) is based on a first-order
optimality criterion for problem~(\ref{subproblem}).  A stronger
assumption on the subproblem solution is made in the following
theorem. Namely, it is assumed that, instead of
requesting~(\ref{modelbaja}) and~(\ref{paradasub}), a global solution
to subproblem~(\ref{subproblem}) is computed. This assumption could be
rather mild in the case that all the subproblems are chosen to be
small dimensional. In this case, it is possible to prove that, in the
limit, suitable $p$th-order optimality conditions are satisfied.
Observe that partial derivatives that are not necessary for computing
Taylor approximations are not assumed to exist at all, let alone to be
continuous.

\begin{teo} \label{teo5}
Suppose that Assumption~\ref{a5} holds and the sequence $\{x^k\}$ is
generated by Algorithm~\ref{cda}.1. Suppose that, at iteration~$k$,
the function $\f$ has as variables $x_i$ with $i \in I_k$, $\Ome$ is
the box $\Omega$ restricted to the variables $i \in I_k$, $\Mbar(x)$
is chosen as the $p$th-order Taylor polynomial of $\f$ defined
in~(\ref{taylorp}), the derivatives involved in (\ref{taylorp}) exist
and are continuous for all $x \in \Omega$, and
Algorithm~\ref{preliminar}.1 computes~$x^+$ as a global minimizer
of~(\ref{subproblem}). Let $K$ be an infinite set of indices such that
$I = I_k$ for all $k \in K$. Let $x^*$ be a limit point of the
sequence $\{x^k\}_{k \in K}$. Then, for all $j \leq p$, $x^*$ is
$j$th-order stationary of problem (\ref{mainprob}) according to
Definition~\ref{def31} and it is also $j$th-order $\sigma$-stationary
for some $\sigma \leq \tau_2 (L+\alpha)$ according to
Definition~\ref{def32} of problem~(\ref{mainprob}).
\end{teo}

\begin{pro}
Consider the problem
\begin{equation} \label{subprolim}
\mbox{Minimize } T_p(x^*, x) + \sigma \|x - x^*\|^{p+1}
\mbox{ subject to } x \in \Omega \mbox{ and } x_i = x_i^* \mbox{ for all } i \notin I.
\end{equation}
By the hypothesis, for all $k \in K$, $x^+$ is obtained as a global
minimizer of
\begin{equation} \label{subprol}
\mbox{Minimize } T_p(x^k, x) + \sigma \|x - x^k\|^{p+1}
\mbox{ subject to } x \in \Omega \mbox{ and } x_i = x_i^k \mbox{ for all } i \notin I,
\end{equation}
for some $\sigma > 0$. Then, by Theorem~\ref{preliminar}.1, $x^{k+1}$
is a global minimizer of~(\ref{subprol}) with $\sigma = \sigma_k \leq
\tau_2 (L + \alpha)$. By (\ref{limxk}), $\lim_{k \in K} x^{k+1} =
\lim_{k \in K} x^k = x^*$. Taking a convenient subsequence, assume,
without loss of generality, that $\lim_{k \in K} \sigma_k = \sigma_*
\leq \tau_2 (L + \alpha)$.  Let $x \in \Omega$ be such that $x_i =
x_i^*$ for all $i \notin I$. Let $z^k \in \Omega$ be such that $z^k_i
= x_i$ for all $i \in I$ and $z^k_i = x_i^k$ for all $i \notin
I$. Then, by the definition of $x^{k+1}$, for all $k \in K$,
\begin{equation} \label{esmenor}
T_p(x^k, x^{k+1}) + \sigma_k \|x^{k+1} - x^k\|^{p+1} \leq
T_p(x^k,z^k) + \sigma_k \|z^k - x^k\|^{p+1}.
\end{equation}
Taking limits for $k \in K$, by the definition of $z^k$, we have that
\begin{equation} \label{latesi}
T_p(x^*, x^*) + \sigma_* \|x^{*} - x^*\|^{p+1} \leq
T_p(x^*, x) + \sigma_* \|x - x^{*}\|^{p+1}.
\end{equation} 
Since $x$ was arbitrary, this implies that $x^*$ is a global solution
of (\ref{subprolim}). Consequently, $x^*$ is also a local solution of
(\ref{subprolim}). Since the Taylor polynomial of order $p$ of
$T_p(x^*, x) + \sigma_* \|x - x^*\|^{p+1}$ coincides with the Taylor
polynomial of order $p$ of $\f$, the thesis is proved.
\end{pro}\\

\noindent
\textbf{Remark~1.} Theorem~\ref{teo5} shows that the convergence of
our CD method is related to an inclusive optimality condition, which
is stronger than every possible exclusive optimality condition.\\
  
\noindent
\textbf{Remark~2.} Note that the hypothesis of Theorem~\ref{teo5}
implies a stronger thesis than the one stated. In fact, we proved
that, in the limit, each partial Taylor polynomial has a global
minimizer. This is interesting because that fact is not a necessary
optimality condition, as it has been shown in the counterexample
exhibited in Section~\ref{background}. However, since \textbf{C3}
implies \textbf{C4} and \textbf{C5}, it turns out that $x^*$ certainly
satisfies the inclusive optimality condition \textbf{C5} according to
Definition~\ref{def33}.

\begin{coro} \label{coroteo5}
Consider the assumptions of Theorem~\ref{teo5} and assume that, for
all $k$,
\[
I_k = \{ \mathrm{mod}(k,n) + 1\}.
\]
If $x^*$ is a limit point of the sequence generated by
Algorithm~\ref{cda}.1, then for all $i = 1, \dots, n$, $x^*_i$ is a
$j$th-order stationary point of the problem
\begin{equation} \label{unaporvez}
\mbox{Minimize } f(x^*_1, \dots, x^*_{i-1}, x_i, x^*_{i+1}, \dots, x^*_n)
\mbox{ subject to } \ell_i \leq x_i \leq u_i
\end{equation}
for all $j \leq p$.
\end{coro}

\begin{pro}
The proof is a direct application of Theorem~\ref{teo5}.
\end{pro}

\section{Complexity} \label{complexity} 

Given a tolerance $\varepsilon > 0$, we wish to know the worst
possible computer effort that we need to obtain an iterate $x$ at
which the objective function is smaller than a given target or the
projected gradient norm $\|g_P(x)\|$ is smaller than $\varepsilon$.
We show that the number of iterations that are needed to obtain
$|[g_P(x^{k+1})]_i| \leq \varepsilon$ for all $i \in I_k$ is, at most,
a constant times $\varepsilon^{-(p+1)/p}$ as in typical high-order
methods. However, obtaining $|[g_P(x^{k+1})]_i| \leq \varepsilon$ for
all $i \notin I_k$ is harder as, for this purpose, we need that
consecutive iterations be close enough. This difficulty is intrinsic
to coordinate descent methods. Powell's example of non-convergence of
CD methods~\cite{powell} satisfies the requirement $|[g_P(x^{k+1})]_i|
\leq \varepsilon$ for all $i \in I_k$ at every iteration but never
satisfies $|[g_P(x^{k+1})]_i| \leq \varepsilon$ for $i \notin
I_k$. Our method converges even in Powell's example because we require
sufficient descent based on regularization but it is affected by
Powell's effect because the number of iterations at which the distance
between consecutive iterates is bigger than a fixed distance grows
with the order~$p$. Then, it is not surprising that our worst-case
complexity bound is significantly worse than
$O(\varepsilon^{-(p+1)/p})$. These results are rigorously proved in
this section and discussed in Section~\ref{discussion}.

\begin{teo} \label{teocomple3}  
Suppose that Assumption~\ref{a5} holds. Let $\ftarget \leq f(x^0)$
and $\varepsilon>0$ be given. Then, the quantity of iterations~$k$
such that
\begin{description}
\item[(i)] $f(x^{k+1}) > \ftarget$ and
\item[(ii)] $|[g_P(x^{k+1})]_i| > \varepsilon$ for some $i \in I_k$
\end{description}
is bounded by
\begin{equation} \label{cota3}
\frac{f(x^0)-\ftarget}{c \, \varepsilon^{(p+1)/p}}, 
\end{equation}
where $c$ only depends on $\alpha$, $\tau_2$, $L$, $p$, and $\theta$.
\end{teo}

\begin{pro}  
By (\ref{cuantoba}) in Theorem~\ref{teo1}, 
\[
f(x^{k+1}) \leq f(x^k) - c \| g_{P,I_k} (x^{k+1})\|^{(p+1)/p},
\]
where $c = \left( \alpha / ( L + \tau_2 (L + \alpha) (p+1) + \theta)
\right)^{(p+1)/p}$. Therefore, if $i \in I_k$,
\[
f(x^{k+1}) \leq f(x^k) - c \left| [g_{P} (x^{k+1})]_i \right|^{(p+1)/p}.
\]   
So, if $\left| [g_P(x^{k+1})]_i \right| > \varepsilon$,
\begin{equation} \label{ocurre2}
f(x^{k+1}) \leq f(x^k) - c \varepsilon^{(p+1)/p}.
\end{equation}
Since the sequence $\{f(x^k)\}$ decreases monotonically, the number of
iterations at which~(\ref{ocurre2}) occurs together with $f(x^{k+1}) >
\ftarget$ cannot exceed $( f(x^0)-\ftarget ) / ( c
\varepsilon^{(p+1)/p} )$. This completes the proof.
\end{pro}

\begin{teo} \label{compledelta}  
Suppose that Assumption~\ref{a5} holds. Let $\ftarget \leq f(x^0)$ and
$\delta > 0$ be given. Then, the quantity of iterations~$k$ such that
$f(x^k) > \ftarget$ and $\|x^{k+1}-x^k\| > \delta$ is bounded by
\begin{equation} \label{cota1}
\frac{f(x^0)-\ftarget}{\alpha \, \delta^{p+1}}.
\end{equation}
\end{teo}

\begin{pro}  
The proof follows directly from (\ref{armijo5}) in Theorem~\ref{teo2}.
\end{pro}

\begin{teo} \label{otro}  
Suppose that Assumption~\ref{a5} holds. Let $\ftarget \leq f(x^0)$,
$\varepsilon > 0$, and $\delta > 0$ be given.  Then, the quantity of
iterations~$k$ such that
\begin{description}
\item[(i)] $f(x^{k+1}) > \ftarget$ and
\item[(ii)] $\|x^{k+1}-x^k\| > \delta$ or $|[g_P(x^{k+1})]_i| >
  \varepsilon$ for some $i \in I_k$
\end{description}
is bounded by
\begin{equation} \label{cota1otra}
  \frac{f(x^0)-\ftarget}{c \, \varepsilon^{(p+1)/p}} +
  \frac{f(x^0)-\ftarget}{\alpha \, \delta^{p+1}},
\end{equation}
where $c$ only depends on $\alpha$, $\tau_2$, $L$, $p$, and $\theta$.
\end{teo}

\begin{pro}  
The proof follows directly from Theorems~\ref{teocomple3}
and~\ref{compledelta}.
\end{pro}\\

We now divide the iterations of Algorithm~\ref{cda}.1 in
\textit{cycles}. Each cycle is composed by~$\mbar$ iterations,
where~$\mbar$ is the one assumed to exist in
Assumption~\ref{aemebar}. Therefore, the successive cycles start at
iterations $x^0, x^{\mbar}, x^{2 \mbar}, \dots, x^{\ell \mbar}, \dots$
The iterates $x^{\ell \mbar+1}, \dots, x^{\ell \mbar+ \mbar}$ are said
to be \textit{produced} at cycle~$\ell$. Iterations $k = \ell \mbar,
\dots, \ell \mbar + \mbar-1$, at which these iterates were produced,
are said to be \textit{internal} iterations of cycle~$\ell$. Each
iteration~$k$ is associated with a set of indices~$I_k$. Due to
Assumption~\ref{aemebar}, for every coordinate $i=1,\dots,n$ and every
cycle~$\ell \geq 0$, there is at least an iteration~$k$ internal to
cycle~$\ell$ such that $i \in I_k$. In other words, all coordinates
are considered in at least an iteration of every cycle. With the
notion of cycle at hand, we can now restate Theorems~\ref{teocomple3},
\ref{compledelta}, and~\ref{otro} as follows.

\begin{teo} \label{compledomingo}  
Suppose that Assumptions~\ref{a5} and~\ref{aemebar} hold. Let
$\ftarget \leq f(x^0)$ and $\varepsilon > 0$ be given. Then, the
quantity of cycles~$\ell$ that contain an internal iteration~$k$ such
that
\begin{description}
  \item[(i)] $f(x^{k+1}) > \ftarget$ and
  \item[(ii)] $|[g_P(x^{k+1})]_i| >  \varepsilon$ for some $i \in I_k$
\end{description}
is not bigger than
\begin{equation} \label{elemenor}
\frac{f(x^0)-\ftarget}{c \, \varepsilon^{(p+1)/p}},  
\end{equation}
where $c$ only depends on $\alpha$, $\tau_2$, $L$, $p$, and $\theta$.
\end{teo}

\begin{pro}
Let $\ell$ be a cycle that contains an internal iteration~$k$
satisfying~(i) and~(ii). By Theorem~\ref{teocomple3}, the quantity of
this type of iteration is bounded by~(\ref{elemenor}); and so the same
bound applies to the quantity of cycles containing an iteration with
these properties. This completes the proof.
\end{pro}

\begin{teo} \label{compledelta2}  
Suppose that Assumptions~\ref{a5} and~\ref{aemebar} hold. Let
$\ftarget \leq f(x^0)$ and $\varepsilon > 0$ be given. Then, the
quantity of cycles~$\ell$ that contain an internal iteration~$k$ such
that $f(x^k) > \ftarget$ and $\|x^{k+1}-x^k\| > \delta$ is bounded by
\begin{equation} \label{cota1repetida}
\frac{f(x^0)-\ftarget}{\alpha \, \delta^{p+1}}.
\end{equation}
\end{teo}

\begin{pro}  
Let $\ell$ be a cycle that contains an internal iteration~$k$ such
that $f(x^k) > \ftarget$ and $\|x^{k+1}-x^k\| > \delta$. By
Theorem~\ref{compledelta}, the quantity of this type of iteration is
bounded by~(\ref{cota1repetida}); and so the same bound applies to the
quantity of cycles containing an iteration with these
properties. This completes the proof.
\end{pro}

\begin{teo} \label{lunes}  
Suppose that Assumptions~\ref{a5} and~\ref{aemebar} hold. Let
$\ftarget \leq f(x^0)$, $\varepsilon > 0$, and $\delta>0$ be
given. Then, the quantity of cycles~$\ell$ that contain an internal
iteration~$k$ such that
\begin{description}
\item[(i)] $f(x^{k+1}) > \ftarget$ and
\item[(ii)] $\|x^{k+1}-x^k\| > \delta$ or $|[g_P(x^{k+1})]_i| >
  \varepsilon$ for some $i \in I_k$
\end{description}
is bounded by
\begin{equation} \label{elemesum}
  \frac{f(x^0)-\ftarget}{c \, \varepsilon^{(p+1)/p}} +
  \frac{f(x^0)-\ftarget}{\alpha \, \delta^{p+1}},
\end{equation}
where $c$ only depends on $\alpha$, $\tau_2$, $L$, $p$, and $\theta$.
\end{teo}

\begin{pro}  
The proof follows directly from Theorems~\ref{compledomingo} and
\ref{compledelta2}.
\end{pro}\\      

The following assumption guarantees that small increments cause small
differences on the projected gradients.

\begin{assump} \label{alips}
There exists $L_g > 0$ such that for all $i=1,\dots,n$ and $x, z \in
\Omega$,
\begin{equation} \label{lipsind}
\left| [  g_P(x)]_i - [g_P(z)]_i \right| \leq L_g \|x - z\|.
\end{equation}
\end{assump}

By the non-expansiveness property of projections,
Assumption~\ref{alips} is satisfied if the gradient of~$f$ satisfies a
Lipschitz condition with constant~$L_g$.

With the tools given by Assumption~\ref{alips} and
Theorem~\ref{lunes}, we are now able to establish a bound on the
number of cycles at which the whole projected gradient is bigger than
a given tolerance.
 
\begin{teo} \label{teocomple5}  
Suppose that Assumptions~\ref{a5}, \ref{aemebar}, and~\ref{alips}
hold. Let $\ftarget \leq f(x^0)$, $\varepsilon > 0$, and $\delta > 0$
be given. Then, there exists a cycle $\ell$, with~$\ell$
exceeding~(\ref{elemesum}) by one in the worst case, such that either
\begin{description}
  \item[(i)] for some iteration~$k$ internal to cycle~$\ell$, we have
    that $f(x^k) \leq \ftarget$ or
  \item[(ii)] for all the iterations $k$ internal to cycle~$\ell$ we have that
    \begin{equation} 
      \left| [g_P(x^{k+1})]_i \right| \leq \varepsilon + \bar m L_g
      \delta \mbox{ for all } i = 1, \dots, n.
    \end{equation}
\end{description}
\end{teo}

\begin{pro}
By Theorem~\ref{lunes}, there exists a cycle $\ell$ that does not
exceeds~(\ref{elemesum}) by more than one such that, for each
iterations~$k$ internal to cycle~$\ell$, either $f(x^{k+1}) \leq
\ftarget$ or
\begin{equation} \label{eso}
\|x^{k+1}-x^k\| \leq \delta \mbox{ and } |[g_P(x^{k+1})]_i| \leq
\varepsilon \mbox{ for all } i \in I_k.
\end{equation}
If there exists an iteration~$k$ internal to cycle~$\ell$ such that
$f(x^{k+1}) \leq \ftarget$, then we are done. So, we assume that, for
all iterations~$k$ internal to cycle~$\ell$, (\ref{eso}) holds. Let $i
\in \{1,\dots,n\}$ be arbitrary. Assumption~\ref{aemebar} implies
that there is an iteration~$k$ internal to cycle~$\ell$ such that $i
\in I_k$ and, thus, by~(\ref{eso}), $|[g_P(x^{k+1})]_i| \leq
\varepsilon$. For any other iterate~$z$ produced at cycle~$\ell$, by
Assumption~\ref{alips}, the triangle inequality, and the first
inequality in~(\ref{eso}), we have that
\[
| [g_P(z)]_i - [g_P(x^{k+1})]_i | \leq L_g \| z - x^{k+1} \| \leq \bar m L_g \delta.
\]
Thus,
\[
| [g_P(z)]_i | \leq \varepsilon + \bar m L_g \delta,
\]
as we wanted to prove.
\end{pro}

\begin{teo} \label{teocomple6}  
Suppose that Assumptions~\ref{a5}, \ref{aemebar}, and~\ref{alips}
hold. Let $\ftarget \leq f(x^0)$, $\varepsilon > 0$, and $\delta > 0$
be given. Then, there exists a cycle $\ell$ of index not larger than
\begin{equation} \label{kepsdelta2}
\frac{f(x^0)-\ftarget}{c \,
  (\varepsilon/2)^{(p+1)/p}} + \frac{f(x^0)-\ftarget}{\alpha \,
  (\varepsilon/(2 \bar m L_g)^{p+1})} + 1,
\end{equation}
where $c$ only depends on $\alpha$, $\tau_2$, $L$, $p$, and $\theta$,
such that, in its first internal iteration~$k$, either $f(x^k) \leq
\ftarget$ or
\begin{equation} 
\left| [g_P(x^{k+1})]_i \right| \leq \varepsilon \mbox{ for all } i =
1, \dots, n.
\end{equation}
\end{teo}

\begin{pro}
The proof follows from Theorem~\ref{teocomple5} replacing
$\varepsilon$ with $\varepsilon/2$ and defining $\delta = \varepsilon
/(2 \bar m L_g)$. Note that the thesis holds for the first iteration
of the cycle because, in fact, due to Theorem~\ref{teocomple5}, it
holds for all its iterations.
\end{pro}\\

The impact of $\bar m$ on the complexity limit is expressed in
formula~(\ref{kepsdelta2}). Note that the second term
of~(\ref{kepsdelta2}) grows proportionally to $\bar m^{p+1}$. If $n$
increases and the size of the subproblems remains bounded, then $\bar
m$ grows proportionately to $n$. Under these conditions, an increase
in the number of iterations proportional to $n^{p+1}$ is
expected. Theorems~\ref{teocomple3}--\ref{teocomple6} give upper
bounds on the number of iterations of Algorithm~\ref{cda}.1. (Bounds
on the number of cycles translate into bounds on the number of
iterations if multiplied by~$\bar m$.) The first term of the sequence
of regularization parameters used in Algorithm~\ref{preliminar}.1 is
$0$. If the corresponding trial point is rejected, the second term is
$\sigma_{\min}$. Then, each time that $\sigma$ needs to be increased,
it is multiplied by a number larger than or equal to
$\tau_1$. Therefore, by definition, the sequence of $\sigma$'s
generated by Algorithm~\ref{preliminar}.1 is bounded from below by the
sequence $0, \, \tau_1^0 \sigma_{\min}, \, \tau_1^1 \sigma_{\min}, \,
\tau_1^2 \sigma_{\min}, \, \tau_1^3 \sigma_{\min}, \, \dots$ Thus, by
Theorem~\ref{teo1}, the number of functional evaluations per call to
Algorithm~\ref{preliminar}.1 at Step~2 of Algorithm~\ref{cda}.1 is
bounded by
\[
\log_{\tau_1}((L+\alpha)/\sigma_{\min}) + 2.
\]
This establishes analogous bounds on the number of functional
evaluations of Algorithm~\ref{cda}.1.

\section{Discussion} \label{discussion}

Theorems~\ref{compledomingo} and~\ref{compledelta2} are complementary
for showing that, eventually, Algorithm~\ref{cda}.1 computes an
iterate $x^k$ such that $\|\g_P(x^k)\|$ is smaller than a given
tolerance; and that this task employs an amount of computer time that
depends on tolerances and problem parameters. In
Theorem~\ref{compledomingo}, we proved that within
$O(\varepsilon^{-(p+1)/p})$ iterations Algorithm~\ref{cda}.1 computes
a sequence (cycle) of $\mbar$ iterates such that, for each $i = 1,
\dots, n$, there is at least one $k$ such that $|[g_P(x^{k+1})]_i|
\leq \varepsilon$. The number of required iterations for this purpose
decreases with~$p$ and tends to $O(1/\varepsilon)$ when~$p$ tends to
infinity. However, this result does not guarantee that the projected
gradient norm is smaller than~$\varepsilon$ at a single iterate. For
this purpose, we need the different iterates within a cycle to be
clustered in a ball of small size. Unfortunately, in order to
guarantee that this happens with tolerance $\delta$, we need,
according to Theorem~\ref{compledelta2}, $O(1/\delta^{p+1})$
iterations. This quantity increases with $p$, which seems to indicate
that, in the worst case, high-order coordinate descent is less
efficient than low-order coordinate descent.

Examples given by Powell in \cite{powell} indicate that, in fact, this
may be the case. In these examples, if coordinate descent is employed
with exact coordinate minimization and cyclic coordinate descent, the
generated sequence has more than one limit point. So, the distance
between consecutive iterations does not tend to zero. This behavior is
not observed if Algorithm~\ref{cda}.1 is applied because the descent
condition~(\ref{armijo4}) implies that $\lim \|x^{k+1}-x^k\| =
0$. However, exact minimization at each iteration evokes the case $p =
\infty$ of Algorithm~\ref{cda}.1 in the sense that the trial point
computed as an exact minimizer satisfies the conditions for accepting
the trial steps for any $p$. So, the conjecture arises that if one
applies Algorithm~\ref{cda}.1 to Powell's examples with different
values of $p$, the resulting sequence, although convergent to a
solution, stays an increasing number of iterations oscillating around
Powell's limiting cycle.

This conjecture is not easy to verify because, except one, Powell's
examples are unstable in the sense that small perturbations cause
convergence to the true minimizers far from the limit spurious cycle.
In any case, we can emulate the application of Algorithm~\ref{cda}.1
to the most famous of Powell's examples (slightly modified here):
\begin{equation} \label{powellex}
\mbox{ Minimize } f(x_1, x_2, x_3) \equiv
   - (x_1 x_2 + x_1 x_3 + x_2 x_3)  +  \sum_{i=1}^3 (|x_i| - 0.1)_+^2. 
\end{equation}
If coordinate descent method employing exact coordinate minimization
and cyclic coordinate descent is applied to problem~(\ref{powellex})
starting from
\[ 
x^0 = (-0.1-\epsilon, 0.1 + \epsilon/2, -0.1-\epsilon/4),
\]
it generates, after six iterations, an iterate~$x^6$ that
corresponds to~$x^0$ with~$\epsilon$ substituted with $\epsilon/64$,
i.e.\
\[
x^6 = (-0.1-\epsilon/64, 0.1 + \epsilon/128, -0.1-\epsilon/256);
\]
and, in general, for all~$k$,
\[
x^{6k} = (-0.1-\epsilon/64^k, 0.1 + \epsilon/(2 \times 64^k),
-0.1-\epsilon/(4 \times 64^k)).
\]
In the intermediate iterations, that are not multiples of~$6$, one has
that
\[
x^{6k+j} = (\pm 0.1 \pm \epsilon/\nu_{k,j}, \pm 0.1 \pm \epsilon/
\times \nu_{k,j}, \pm 0.1 \pm \epsilon/\nu_{k,j})
\]    
where $\nu_{k,j} \leq 4 \times 64^{k+1}$ for all $k, j$.

Now, we wish to show that this sequence could be generated by
Algorithm~~\ref{cda}.1. Moreover, for any given $p$, we wish to know
how many iterations are necessary to obtain consecutive iterations
such that $\|x^{k+1}-x^k\| \leq 0.01$. Let
\[ 
x^0 = (-0.1-\epsilon, 0.1 + \epsilon/2, -0.1-\epsilon/4).
\]
The global minimizer of $f(x_1, x_2, x_3)$ subject to $x_2 = x_2^0$
and $x_3 = x_3^0$ is
\[
z^0 = (0.1 + \epsilon/8, 0.1 + \epsilon/2, -0.1-\epsilon/4).
\]
(The iterate $x^1$ in the Powell's sequence is given by $x^1=z^0$, but
we preserve the notation $z^0$ for the sake of simplicity.) On the one
hand,
\[
f(x^0) = - (x_1^0 x_2^0 + x_1^0 x_3^0 + x_2^0 x_3^0)  +  \sum_{i=1}^3 (|x_i^0| - 0.1)_+^2. 
\]
On the other hand, since $z^0_2 = x^0_2$ and $z^0_3 = x^0_3$,
 \[
f(z^0) = - (z_1^0 x_2^0 + z_1^0 x_3^0 + x_2^0 x_3^0) + (|z_1^0| -
0.1)_+^2 + \sum_{i=2}^3 (|x_i^0| - 0.1)_+^2.
\]
Therefore,
\[
f(x^0) - f(z^0) = (z_1^0 - x_1^0) (x^2_0 + x^3_0) + (|x_1^0| -
0.1)_+^2 - (|z_1^0| - 0.1)_+^2.
\]
Thus,
\[
\resizebox{\textwidth}{!}{$\begin{array}{rcl}
f(x^0) - f(z^0) &=& ((0.1+\epsilon/8) - (-0.1-\epsilon)) (\epsilon/2 -
\epsilon/4) + (|-0.1 - \epsilon| - 0.1)_+^2 - (|0.1 + \epsilon/8| -
0.1)_+^2\\[2mm]
&=& (0.2+ 9 \epsilon/8) \epsilon/4 + \epsilon^2 - \epsilon^2/64 = 0.2
\epsilon/4 + 9 \epsilon^2/32 + \epsilon^2 - \epsilon^2/64\\[2mm]
&=& 0.2 \epsilon/4 + 9 \epsilon^2/32 + \epsilon^2 - \epsilon^2/64 = 0.2
\epsilon/4 + 81 \epsilon^2/64 \geq \epsilon/20.
\end{array}$}
\]

Consider Algorithm~\ref{cda}.1 using $f(x)$ as the model of the
objective function. We must verify whether (\ref{modelbaja}),
(\ref{paradasub}), and~(\ref{armijo4}) are satisfied with $\xtrial =
z^0$. Trivially, for $\sigma = 0$, (\ref{modelbaja})
and~(\ref{paradasub}) hold by the definition of the model and the fact
that $z^0$ is a global minimizer. In order to show
that~(\ref{armijo4}) also holds, let as assume that $\epsilon < 0.1$
and $2^{p+1} \geq 20 \alpha/\epsilon$, i.e.\ $\alpha / 2^{p+1} \leq
\epsilon/20$. So, by the calculations above,
\[
f(x^0) - f(\xtrial)  \geq   \alpha /2^{p+1}.
\]   
Since $\epsilon < 0.1$, we have that $\|\xtrial - x^0\| \leq
0.5$. Thus,
\[
f(x^0) - f(\xtrial)  \geq   \alpha  \|\xtrial - x^0\| ^{p+1}.
\]   
This implies (\ref{armijo4}). Therefore, a sufficient condition for
the acceptance of~$x^1=z^0$ as an iterate of Algorithm~\ref{cda}.1 is
\[
\alpha / 2^{p+1} \leq \frac{\epsilon}{20 \times 4 \times 64^k}.
\]   
In other words,
\[
20 \times 4 \times 64^k   \alpha \leq   \epsilon  2^{p+1}.
\]  
Taking logarithms, this condition is
\[
\log_2 80 + 6 k + \log_2 \alpha \leq p+1.
\]
That is, if
\[
k_0  \leq ( p+1 - \log_2 80 - \log_2 \alpha ) / 6,
\]   
the first $k_0$ iterations of Algorithm~\ref{cda}.1 will reproduce the
cycling example of Powell. In all these iterations we have that
$\|x^{k+1}-x^k\| \geq 0.1$. Note that~$k_0$ tends to infinity as~$p$
tends to infinite, as we wanted to show. In addition, note also
that~$k_0$ tends to infinity as~$\alpha$ tends to zero, which reflects
the obvious fact that, if we are more tolerant with the acceptance of
the trial point, the probability of staying around Powell's six-points
cycle increases.
                            
It is not sensible to decide about usefulness of algorithms based only
on theoretical convergence or complexity results. Since these results
deal with worst-case behavior the possibility exists that a class of
problems in which practitioners are interested always exhibit
characteristics that exclude extreme unfortunate cases. However, it is
pertinent to examine pure mathematical properties in order to foster
unexpected good or bad computer behaviors.

\begin{enumerate}

\item Many optimization users believe that if a smooth function has a
  minimizer at a point~$x^*$, then this point is a local minimizer of
  all its Taylor polynomials.  This is true only if the dimension $n$
  is equal to $1$. For arbitrary $n$, it is true only up to second
  order polynomials. Examples that illustrates this phenomenon have
  been given in this paper with the purpose of justifying adequate
  high-order optimality conditions (for example, $f(x_1, x_2) = x_2^2
  - x_1^2 x_2 + x_1^4$). This fact implies that, in the vicinity of a
  global minimizer, a high-order algorithm may try to find
  improvements far from the current point, being subject to a painful
  sequence of ``backtrackings" before obtaining descent. Does this
  imply that only quadratic approximations are useful in the
  minimization context?  It is too soon to give a definite response to
  this question.

\item Our regularization approach for CD-algorithms makes it
  impossible to exhibit the cyclic behavior of Powell's examples
  \cite{powell}. The reason is that, under regularization descent
  algorithms, the difference between consecutive iterates tends to
  zero.  However, it seems to be possible that convergence to zero of
  consecutive iterates could be very slow, as predicted by complexity
  results. Is this an argument for discarding high-order CD
  algorithms? We believe that the answer is no, as far as the use of
  CD algorithms is, in general, motivated by the structure of the
  problems, which in some sense should evoke some degree of
  separability.  Moreover, since high-order models are also low-order
  models one can use high-order associated with a small $p$ in
  (\ref{modelbaja}), (\ref{paradasub}), and (\ref{armijo4}). In other
  words, if $1 \leq q < p$, then the conditions that define a model of
  order~$q$ are satisfied by models of order $p$. Therefore, we may
  use models of order~$p$ associated with the regularization required
  by models of order~$q$. For example, we may use a second-order model
  associated to quadratic regularization preserving first-order
  convergence results and the corresponding complexity.

\item It is interesting to consider the case in which we use $f(x)$ as
  a model for $f(x)$. In this case, high-order analysis makes a lot of
  sense. In fact, efficient algorithms for finding global minimizers
  of functions of one variable exist, a possibility that decreases
  very fast as the number of variables grow. Moreover high-order
  one-dimensional models are certainly affordable and many numerical
  analysis papers handle efficiently the problem of minimizing or
  finding roots of univariate polynomials \cite{petkovic}.  Recall
  that, in this case, the model satisfies the approximation
  requirements for every value of $p$. Therefore we may choose the
  value of $p$ that promises better efficiency, which, according to
  Theorem~\ref{teocomple6}, should be $p=1$ giving complexity
  $O(\varepsilon^{-2})$ as gradient-like methods.

\item In most practical situations one is interested in finding global
  minimizers or, at least, feasible points at which the objective
  function value is smaller than a given~$\ftarget$. Complexity and
  convergence analyses in the nonconvex world concern only the
  approximation to stationary points although every practical
  algorithm must be devised taking into account the global implicit
  goal. It turns out that low coordinate global strategies for finding
  initial points are available in many real-life problems. These
  strategies fit well with CD algorithms as we will illustrate in
  Section~\ref{exp}.

\item The reader will observe that in our experiments we used $p=2$,
  in spite that, according to the complexity results, the optimal $p$
  should be $1$. The reason is that, as we stated in the convergence
  section, the employment of $p=2$ guarantees convergence to points
  that satisfy second order conditions that are not guaranteed by
  $p=1$. Moreover, subproblems with $p=2$ are computationally
  affordable in the applications considered. Summing up, we could say
  that making an informal balance regarding theoretical results, using
  $p=2$ should be the default choice for practical applications.
\end{enumerate}

\section{Implementation and experiments} \label{exp}

This section illustrates with numerical experiments the applicability
of Algorithm~\ref{cda}.1. The Multidimensional Scaling (MS) problem
\cite{coxcox,mead,torgerson} adopted for the experiments is described
in Section~\ref{exp1}. Implementation details of
Algorithms~\ref{preliminar}.1 and~\ref{cda}.1 are described in
Section~\ref{exp2}. Problem-dependent strategies for generating an
initial point and for generating a sequence of improved initial points
are described in Section~\ref{exp3}. The computational results are
shown in Section~\ref{exp4}.

\subsection{Multidimensional Scaling problem} \label{exp1}

Multidimensional Scaling methods emerged as statistical tools in
Psychophysics and sensory analysis.  The MS problem considered in this
section may be described in the following way: Let $x_1,\dots,x_{n_p}
\in \R^d$ be a set of unknown points. Let $D=(d_{ij}) \in \R^{n_p
  \times n_p}$ be such that $d_{ij} = \| x_i - x_j\|$; and assume that
only entries $\{ d_{ij} \;|\; (i,j) \in S \}$ for a given $S \subset
\{ 1,\dots,n_p \} \times \{ 1,\dots,n_p \}$ are known. (Of course, $D$
is symmetric, $d_{ii}=0$, and $(i,j) \in S$ if and only if $(j,i) \in
S$.) Then the MS problem consists of finding $x_1,\dots,x_{n_p}$ such
that $\| x_i - x_j\| = d_{ij}$ for all $(i,j) \in S$. Glunt, Hayden,
and Raydan~\cite{ghr} were the first to apply unconstrained continuous
optimization tools to the nowadays called Molecular Distance Geometry
Problem (MDGP), as defined in \cite{carlile,carlile2} in a
Multidimensional Scaling context.  This problem appears when points
$x_1,\dots,x_{n_p}$ correspond to the positions of atoms in a molecule
and distances not larger than 6 Angstroms (i.e.\ $6 \times 10^{-10}$
meters) are obtained via nuclear magnetic resonance (NMR)
\cite{artacho}. This problem can be modeled as the following
unconstrained nonlinear optimization problem
\begin{equation} \label{MDGP}
\Minimize_{x_1,\dots,x_{n_p} \in \R^d} f(x_1,\dots,x_{n_p})
:= \frac{1}{|S|} \sum_{(i,j) \in S} \left( \| x_i - x_j \|_2^2 - d_{ij}^2 \right)^2.
\end{equation}

\subsection{Implementation details} \label{exp2}

If we wish to apply Algorithm~\ref{cda}.1 to the MDGP problem, it
arises quite naturally to associate at iteration~$k$ the set $I_k$
with the components of a point $x_{\ell(k)} \in \R^d$ for some
$\ell(k)$ between $1$ and $n_p$. Specifically, if we define
$x=(x_1^T,\dots,x_{n_p}^T)^T \in \R^n$ with $n:=d \, n_p$, then at
iteration~$k$ we can define
\begin{equation} \label{Ik}
  I_k = \left\{ (\ell(k)-1)d+1,\dots,(\ell(k)-1)d+d \right\}
  \mbox{ with } \ell(k)=\mathrm{mod}(k,n_p)+1,
\end{equation}
or any alternative choice of $\ell(k) \in \{1,\dots,n_p\}$.
%In the MDGP, the idea of associating the set $I_k$ of iteration~$k$
%with the components of a point $x_k \in \R^d$ that is modified at each
%iteration of Algorithm~\ref{cda}.1 is quite natural, i.e.\ we define
%$x=(x_1^T,\dots,x_{n_p}^T)^T \in \R^n$ with $n:=d \, n_p$ and, at
%iteration $k$,
%\begin{equation} \label{Ik}
%I_k = \left\{ (\ell(k)-1)d+1,\dots,(\ell(k)-1)d+d \right\}
%\mbox{ with } \ell(k)=\mathrm{mod}(k,n_p)+1,
%\end{equation}
%or any alternative choice of $\ell(k) \in \{1,\dots,n_p\}$.
This is equivalent to say that, at iteration $k$, the subproblem
considered at Step~2 of Algorithm~\ref{cda}.1 is given by
\begin{equation} \label{ell}
\Minimize_{z \in \R^d} \underline{f}(z),
\end{equation}
where $\underline{f} : \R^d \to \R$ is defined as
\begin{equation} \label{ell2}
\underline{f}(z) := \frac{1}{|S|} \left[
\sum_{(i,j) \in S \setminus S(\ell(k))} \left( \| x_i - x_j \|_2^2 - d_{ij}^2 \right)^2 + 
2 \sum_{(i,\ell(k)) \in S} \left( \| x_i - z \|_2^2 - d_{i,\ell(k)}^2 \right)^2 \right],
\end{equation}
$S(\ell(k)):= \{ (i,j) \in S \; | \; i=\ell(k) \mbox{ or }
j=\ell(k)\}$, and $\ell(k)$ is given by~(\ref{Ik}). Note that the time
complexity for evaluating~$f$ is $O(d|S|)$; while, since the first
summation in~(\ref{ell2}) does not depend on $z$, the time complexity
for evaluating~$\underline{f}$ is, in average $O(d|S|/n_p)$.

For approximately solving~(\ref{ell}) in Algorithm~\ref{preliminar}.1,
we consider a second-order Taylor expansion of~$\underline{f}$ at
$\bar x = x_{\ell(k)}^k \in \R^d$, i.e.\
\begin{equation} \label{modelellmenosuno}
M_{\bar x}(z) := \underline{f}(\bar x)
+ \nabla \underline{f}(\bar x)^T (z - \bar x)
+ (z - \bar x)^T \nabla^2 \underline{f}(\bar x)^T (z - \bar x).
\end{equation}
This means that the underlying model-based subproblem, when
Algorithm~\ref{preliminar}.1 is used at Step~2 of the $k$th iteration
of Algorithm~\ref{cda}.1 is given by
\begin{equation} \label{modelell}
\Minimize_{z \in \R^d} M_{\bar x}(z) + \sigma \| z - \bar x \|^3.
\end{equation}
Since problem~(\ref{MDGP}) is unconstrained, i.e.\ $\Omega=\R^n$,
subproblems~(\ref{ell}) and model-based subproblems~(\ref{modelell})
are unconstrained as well. Thus, if in~(\ref{modelell}) and, in
consequence, in~(\ref{armijo4}), for $x \in \R^d$, we consider $\| x
\|$ as $\| x \|_3 := ( \sum_{i=1}^d |x_i|^3 )^{1/3}$, then the
\textit{global} minimizer of~(\ref{modelell}) can be easily obtained
at the expense of a single factorization of $\nabla^2
\underline{f}(\bar x) \in \R^{d \times d}$,
see~\cite{bmbunch,bras,mr2014,mr2015}. (When $\sigma=0$,
(\ref{modelell}) may have no solution. This case can be detected with
the same cost as well.)  Since the exact global minimizer $\xtrial$
of~(\ref{modelell}) is being computed at Step~2 of
Algorithm~\ref{preliminar}.1, (\ref{modelbaja}) and~(\ref{paradasub})
always hold, for any $\theta>0$; thus, in the implementation, their
verification can be ignored.

\subsection{Initial guess and multistart strategy} \label{exp3}

As shown in Section~\ref{cda}, Algorithm~\ref{cda}.1 has convergence
properties towards stationary points which, probably, are local
minimizers. Obviously, as we are interested in finding {\it global}
minimizers of MDGP, we need suitable strategies for choosing initial
approximations. We employed the combination of two different
strategies for this purpose. On the one hand, an initial guess
suggested in~\cite{fang} was adopted. On the other hand, we devised a
new coordinate descent procedure based on the structure of MDGP. The
Fang-O'Leary strategy~\cite{fang}, based on shortest paths over an
underlying graph, is a strategy for computing a single initial
solution. Starting from that solution, our new coordinate descent
procedure is used iteratively to make successive improvements on the
Fang-O'Leary initial point. At each improvement,
Algorithm~~\ref{cda}.1 is run to find a local solution.

In order to describe the Fang-O'Leary strategy~\cite{fang}, consider
the weighted graph $G=(\{1,\dots,n_p\},S)$ in which the weight of an
edge $(i,j)$ is given by $d_{ij}$. We assume this graph is
connected. Otherwise, the molecule's structure can not be recovered;
and problem~(\ref{MDGP}) can be decomposed in as many independent
problems as connected components of the graph~$G$ in order to recover
partial structures. Let $\bar S = \{1,\dots,n_p\} \times
\{1,\dots,n_p\} \setminus S$, i.e.\ $\bar S$ corresponds to the
missing arcs in $G$ or, equivalently, the unknown entries of $D$. For
each $(i,j) \in S$, define $\tilde d_{ij}=d_{ij}$; and for each $(i,j)
\in \bar S$, define $\tilde d_{ij}$ as the weight of the shortest path
between $i$ and $j$ in $G$. Matrix $\tilde D = (\tilde d_{ij})$ is a
distance matrix that completes $D$; but with high probability it is
\textit{not} an Euclidean distance matrix. Computing $\tilde D$
requires $O(n_p^2)$ space and has time complexity $O(n_p^3)$ (using
the Floyd-Warshall algorithm as suggested in~\cite{fang}), which can
be an issue for instances with large $n_p$. Obtaining points
$x_1^0,\dots,x_{n_p}^0 \in \R^d$ from $\tilde D$ requires to compute
the $d$ largest positive eigenvalues of the matrix ${\cal T}(\tilde
D)$ given by ${\cal T}(\tilde D):= - \half J \tilde D J$, where
$J:=I-\frac{1}{n} e e^T$ and $e=(1,\dots,1)^T$. If the truncated
spectral decomposition of ${\cal T}(\tilde D)$ is given by $U \Delta_d
U^T$ then the initial point $x^0=((x_1^0)^T,\dots,(x_{n_p}^0)^T)^T$ is
given by $X=(x_1^0,\dots,x_{n_p}^0)=U \Delta_d^{1/2}$. If the matrix
${\cal T}(\tilde D)$ has only $\underline{d} < d$ positive
eigenvalues, then computed points are in $\R^{\underline{d}}$ and
their last $d-\underline{d}$ components can be completed with
zeros. In~\cite{fang}, alternative initial guesses are obtained by
perturbations of matrix~$\tilde D$ and/or by stretching the computed
points $x_1^0,\dots,x_{n_p}^0$.

Our coordinate-descent strategy for choosing the initial approximation
to the solution of MDGP is inspired on the structure of local
solutions. Consider a point $p \in \R^3$ and three other points $q_1,
q_2, q_3 \in \R^3$ such that the distances from $p$ to $q_i$,
$i=1,2,3$, are known, i.e., $(p,q_i) \in S$ for $i=1,2,3$. Assume, in
addition, that the required distances are satisfied, i.e., that
$\|p-q_i\|$ is equal to the corresponding value in matrix $D$ for
$i=1,2,3$. Assume that there is an additional point $q_4$ for which
its known distance $d(p,q_4)$ to $p$ is \textit{not}
satisfied. Assume, in addition, that $(\|r(p) - q_4\|_2^2 -
d(p,q_4)^2)^2 < (\|p - q_4\|_2^2 - d(p,q_4)^2)^2$, where $r(p)$ is the
reflection of $p$ on the plane determined by $q_i$, $i=1,2,3$. If
there were no more points in the problem, replacing $p$ by $r(p)$,
would produce a reduction in the objective function.  Our coordinate
descent algorithm with a coordinate-descent strategy for choosing
initial points is described in Algorithm~\ref{exp}.1. The
coordinate-descent strategy for initial approximations, based on this
intuition, is described at Step~4 of Algorithm~\ref{exp}.1.\\

\noindent
\textbf{Algorithm~\ref{exp}.1.} Assume $\hat x$ is a given arbitrary
initial point (that might be obtained using the Fang-O'Leary technique
described above).
\begin{description}[itemsep=0pt]
\item[\textbf{Step~1.}] Using $\hat x$ as initial guess, run
  Algorithm~\ref{cda}.1 until the obtention of an iterate $\tilde x$
  such that $f(\tilde x) \leq \ftarget$ or such that its projected
  gradient is small enough according to criteria given below.

\item[\textbf{Step~2.}] If $f(\tilde x) \leq \ftarget$ then
  \textbf{stop} declaring that $\tilde x$ is a global minimizer up to
  the precision given by $\ftarget$. Otherwise, update $\hat x$ by
  means of the coordinate-descent strategy in Step~3 below.
\item[\textbf{Step~3.}] For $j=1,\dots,n_p$ execute Steps~3.1--3.2.
\item[\textbf{Step~3.1.}] Let $\hat f_j := \sum_{(i,j) \in S} ( \|
  \hat x_i - \hat x_j \|_2^2 - d_{ij}^2 )^2$.
\item[\textbf{Step~3.2.}] For every triplet $(i_1,i_2,i_3)$ such that
  $(i_1,j), (i_2,j), (i_3,j) \in S$, in an arbitrary order, if
  \[
  \sum_{(i,j) \in S} ( \| \hat x_i - r(\hat x_j) \|_2^2 - d_{ij}^2 )^2 <  \hat f_j,
  \]
  where $r(\hat x_j)$ is the reflection of $\hat x_j$ on the plane
  determined by $\hat x_{i_1}$, $\hat x_{i_2}$, and $\hat x_{i_3}$,
  then update $\hat x_j \leftarrow r(\hat x_j)$. (Note that $\hat f_j$
  is not updated at this point. This means that a sequence of
  reflections can be applied to $\hat x_j$, with a non-monotone
  behavior of $f$, provided it improves the ``reference value'' $\hat
  f_j$.)
\item[\textbf{Step~4.}] If $\hat x$ was not updated at Step~3, then
  stop returning~$\tilde x$. (Note that $\ftarget$ was not reached in
  this case.) Otherwise, go to Step~1.
\end{description}

At Step~1 of Algorithm~\ref{exp}.1, we consider that ``the projected
gradient is small enough'' if, during $n_p$ consecutive iterations of
Algorithm~\ref{cda}.1, we have that ``the final $\sigma$'' of
Algorithm~\ref{preliminar}.1 is larger than $10^{20}$ or $f(x^{k+1})
\not\le f(x^k) - 10^{-8} \min\{1, |f(x^k)|\}$. By (\ref{armijo5}),
(\ref{crucial2}) and the boundedness of~$\sigma$, these are
practical symptoms of stationarity.

\subsection{Computational results} \label{exp4}

We implemented Algorithms~\ref{preliminar}.1, \ref{cda}.1,
and~\ref{exp}.1 in Fortran. In the numerical experiments, we
considered, $\alpha=10^{-8}$, $\sigma_{\min}=10^{-8}$, and
$\tau_1=\tau_2=100$, and $\ftarget = 10^{-10}$. All tests were
conducted on a computer with a 3.5 GHz Intel Core i7 processor and
16GB 1600 MHz DDR3 RAM memory, running macOS High Sierra (version
10.13.6). Code was compiled by the GFortran compiler of GCC (version
8.2.0) with the -O3 optimization directive enabled.

The Research Collaboratory for Structural Bioinformatics (RCSB)
Protein Data Bank~\cite{RCSBPDB} is an open access repository that
provides access to 3D structure data for large biological molecules
(proteins, DNA, and RNA). There are more than $167{,}000$ molecules
available. In~\cite{fang}, where Newton and quasi-Newton methods are
applied to problem~(\ref{MDGP}), six protein molecules are considered,
namely, 2IGG, 1RML, 1AK6, 1A24, 3MSP, and 3EZA (see~\cite[Table~6.9,
  p.20]{fang}); while in~\cite{artacho}, where the Douglas–Rachford
method is applied, other six protein molecules are considered, namely,
1PTQ, 1HOE, 1LFB, 1PHT, 1POA, and 1AX8 (see~\cite[Table~1,
  p.313]{artacho}). In the first work, only protein atoms (identified
with ATOM in the molecule file) were considered; while in the second
work there were considered protein atoms plus atoms in small molecules
(identified with HETATM in the protein molecule file). In the current
work, both options were considered. Following~\cite{fang}, for each
protein molecule, when multiple structures are available, only the
first one was considered. Each molecule is given as the set of 3D
coordinates of its atoms. An instance of problem~(\ref{MDGP}) is built
by computing a complete Euclidean distance matrix and then eliminating
distances larger than 6 Angstroms. Since not all molecules have atoms
in small molecules, we arrived to eighteen different
instances. Table~\ref{tab1} shows, for each instance, the number of
variables $n$ of the optimization problem~(\ref{MDGP}), the number of
atoms~$n_p$, the number of distances considered to be known~$|S|$, and
the CPU time in seconds required to construct the initial guess~$x^0$
using the Fang-O'Leary strategy~\cite{fang}.

\begin{table}[ht!]
\begin{center}
\begin{tabular}{|c|c|cccrr|}
\hline
& & Molecule & $n$ & $n_p$ & \multicolumn{1}{c}{$|S|$} & Time $x^0$ \\  
\hline
\hline
\multirow{18}{*}{\rotatebox{90}{
\begin{tabular}{c}
Points may correspond to protein\\
atoms (ATOM) only or to protein atoms\\
plus atoms in small molecules (HETATM)
\end{tabular}}}
 & \multirow{12}{*}{\rotatebox{90}{ATOM only}}
 & 1ptq &  1,206 &   402 &  14,176 (8.79\%) &   0.21 \\
&& 1hoe &  1,674 &   558 &  20,356 (6.55\%) &   0.49 \\
&& 1lfb &  1,923 &   641 &  22,870 (5.57\%) &   0.70 \\
&& 1pht &  2,433 &   811 &  35,268 (5.37\%) &   1.41 \\
&& 1poa &  2,742 &   914 &  33,966 (4.07\%) &   2.03 \\
&& 2igg &  2,919 &   973 &  62,574 (6.62\%) &   2.54 \\
&& 1ax8 &  3,009 & 1,003 &  37,590 (3.74\%) &   2.76 \\
&& 1rml &  6,192 & 2,064 & 153,660 (3.61\%) &  24.14 \\
&& 1ak6 &  8,214 & 2,738 & 224,568 (3.00\%) &  52.04 \\
&& 1a24 &  8,856 & 2,952 & 212,364 (2.44\%) &  64.90 \\
&& 3msp & 11,940 & 3,980 & 262,876 (1.66\%) & 157.90 \\
&& 3eza & 15,441 & 5,147 & 356,544 (1.35\%) & 335.84 \\
\cline{2-7}
 & \multirow{6}{*}{\rotatebox{90}{{\scriptsize ATOM+HETATM}}}
 & 1ptq &  1,212 &   404 &  14,370 (8.83\%) &   0.21 \\
&& 1hoe &  1,743 &   581 &  21,422 (6.36\%) &   0.55 \\
&& 1pht &  2,964 &   988 &  44,542 (4.57\%) &   2.59 \\
&& 1poa &  3,201 & 1,067 &  41,034 (3.61\%) &   3.23 \\
&& 1ax8 &  3,222 & 1,074 &  40,866 (3.55\%) &   3.29 \\
&& 1rml &  6,273 & 2,091 & 156,550 (3.58\%) &  23.90 \\
\hline
\end{tabular}
\end{center}
\caption{Description of the instances built with the molecules
  considered in~\cite{artacho} or~\cite{fang}.}
\label{tab1}
\end{table}

Note that considered instances are \textit{gedanken} in the sense that
points $\bar x_1, \dots, \bar x_{n_p} \in \R^3$ such that $f(\bar
x)=0$ with $\bar x^T = (\bar x_1^T, \dots, \bar x_{n_p}^T)^T$ are
known. Thus, given $x^*$ such that $f(x^*) \approx 0$, we may wonder
whether $x^*$ is close to $\bar x$. The answer to this question is
``Not necessarily.'' since any rotation or translation of $\bar x$
also annihilates~$f$. So the question would be ``How close is $x^*$ to
$\bar x$ after performing the appropriate rotations and
translations?''. The answer to this question is obtained by solving an
orthogonal Procrustes problem. Let $\bar X = (\bar x_1,\dots,\bar
x_{n_p})$ and $X^* = (x_1^*, \dots, x_{n_p}^*) \in \R^{3 \times
  n_p}$. It is easy to see that matrices $\bar XJ$ and $X^* J$ have
their centroid at the origin, since $\bar XJ e = X^*J e = 0$. (Recall
that $J = I - \frac{1}{n_p} ee^T$ and $e=(1,\dots,1)^T$.) The
orthogonal Procrustes problem consists in finding an orthogonal matrix
$Q \in \R^{3 \times 3}$ which most closely maps $X^* J$ to $\bar X J$,
i.e.\
\[
Q = \argmin_{R \in \R^{3 \times 3}} \| R X^* J - \bar X J \|_F^2 \mbox{ subject to } R R^T = I.
\]
This problem has a closed form solution given by $Q=VU^T$, where $U
\Sigma V^T$ is the singular value decomposition of the matrix $C:=X^*
J (\bar X J)^T$. Thus, the measure we were looking for is given by
\[
E(x^*) := \max_{\{j=1,\dots,n_p\}} \left\{ E(x_j^*) \right\},
\]
where
\begin{equation} \label{errorj}
E(x_j^*) := \frac{\| [Q X^* J - \bar X J]_j \|_{\infty}}{\max\{1,\|
  [\bar X J]_j \|_{\infty}\}},
\end{equation}
and $[A]_j$ denotes the $j$th column of matrix $A$.

\begin{figure}[ht!]
\begin{center}
\begin{tabular}{ccc}
\includegraphics[scale=0.20]{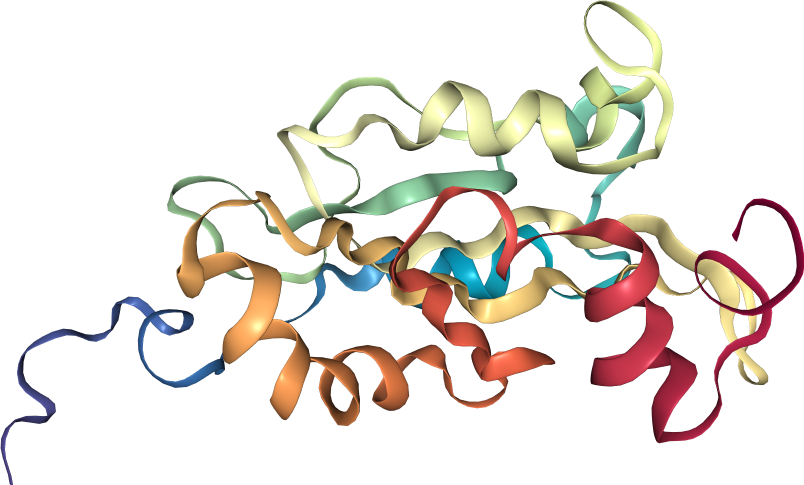} &
\includegraphics[scale=0.14]{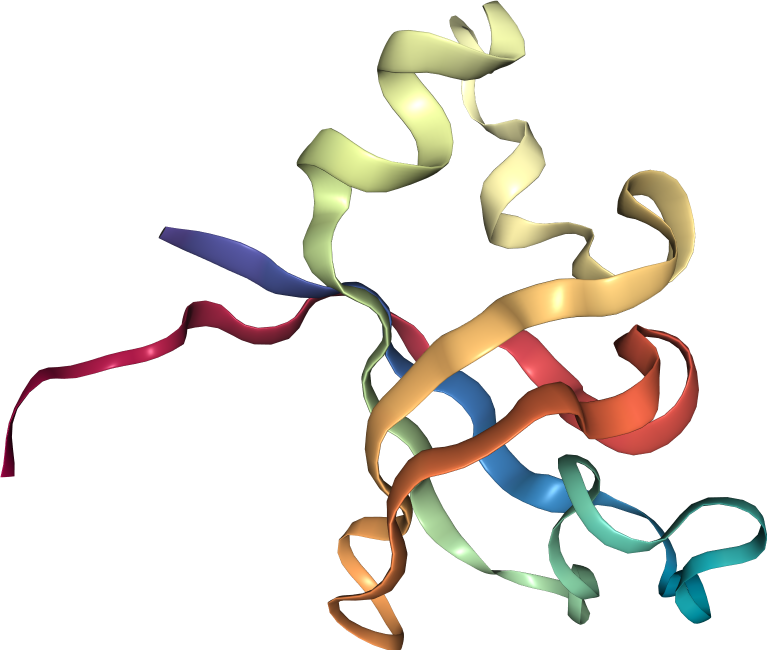} &
\includegraphics[scale=0.14]{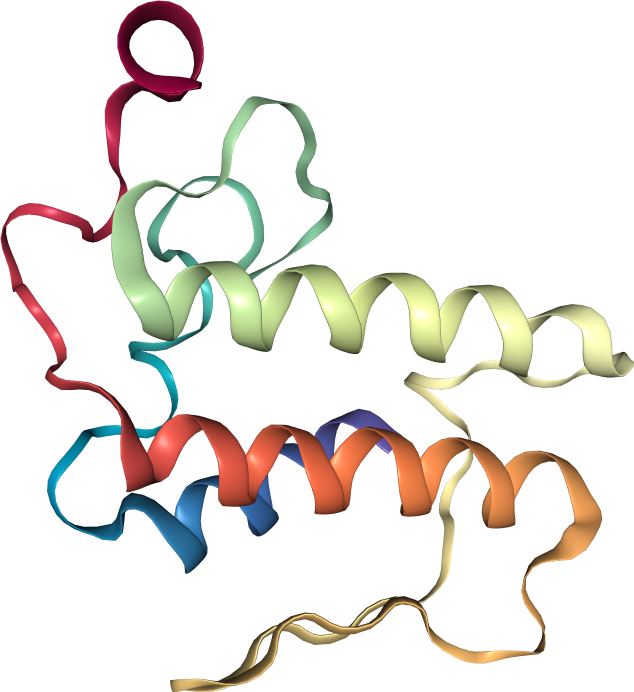} \\
1AK6 & 1PHT & 1POA
\end{tabular}
\end{center}
\caption{Representation of molecules 1AK6, 1PHT, and 1POA for which
  Coordinate Descent found a global minimizer; while SPG and Gencan
  found a local minimizer.}
\label{fig2}
\end{figure}

\begin{figure}[ht!]
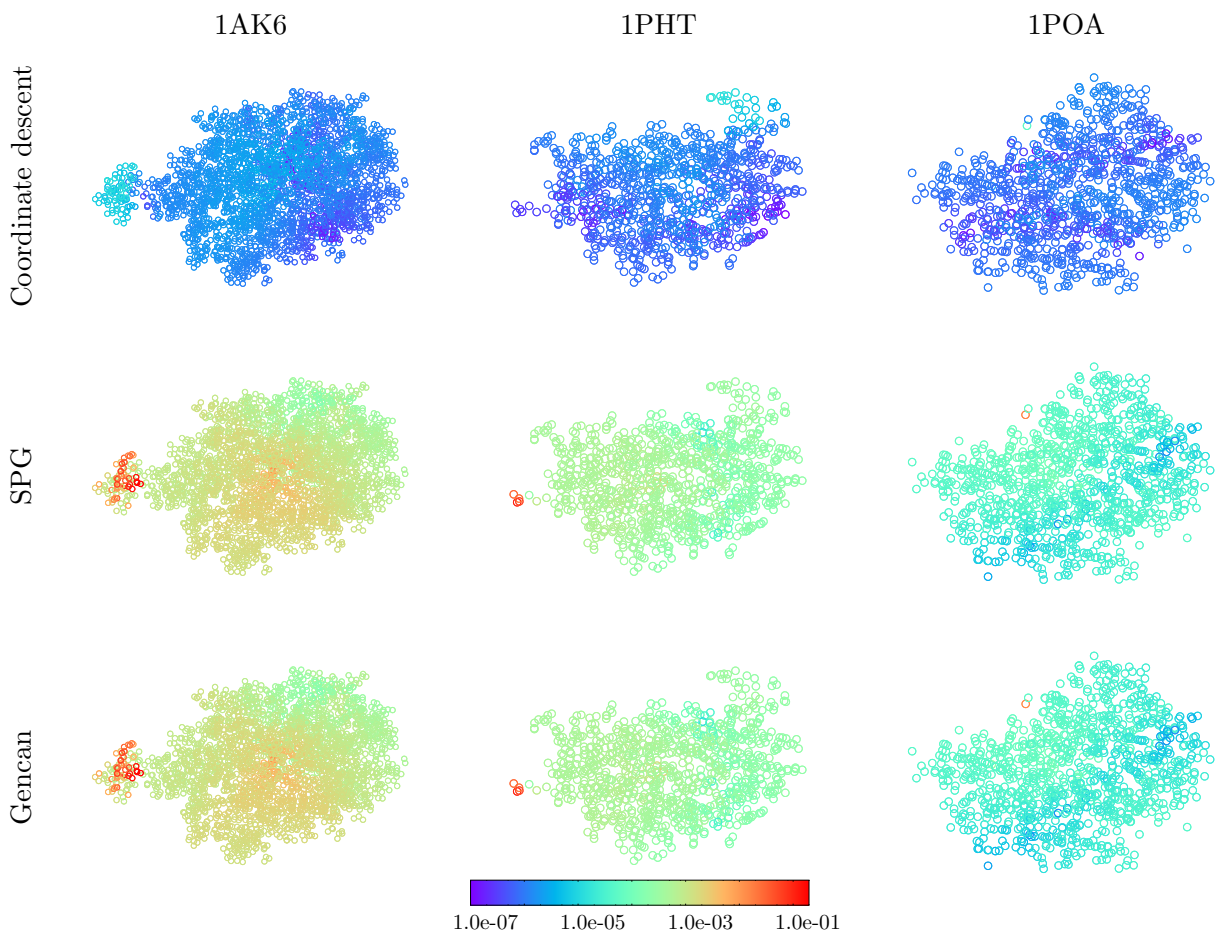

\begin{center}
\begin{tabular}{cccc}
& 1AK6 & 1PHT & 1POA \\[2mm]
\rotatebox{90}{$\phantom{aa}$Coordinate descent} &
\resizebox{5cm}{!}{\input{aabmmfig3aX.tex}} &
\resizebox{5cm}{!}{\input{aabmmfig3bX.tex}} &
\resizebox{5cm}{!}{\input{aabmmfig3cX.tex}} \\
\rotatebox{90}{$\phantom{aaaaaaaa}$SPG} &
\resizebox{5cm}{!}{\input{aabmmfig3dX.tex}} &
\resizebox{5cm}{!}{\input{aabmmfig3eX.tex}} &
\resizebox{5cm}{!}{\input{aabmmfig3fX.tex}} \\
\rotatebox{90}{$\phantom{aaaaaa}$Gencan} &
\resizebox{5cm}{!}{\input{aabmmfig3gX.tex}} &
\resizebox{5cm}{!}{\input{aabmmfig3hX.tex}} &
\resizebox{5cm}{!}{\input{aabmmfig3iX.tex}} \\
\end{tabular}
\end{center}
\caption{Molecules 1AK6, 1PHT, and 1POA for which Coordinate Descent
  found a global minimizer; while SPG and Gencan found a local
  minimizer. To the naked eye, solutions would appear to be
  indistinguishable. Therefore, the figures show, for each point
  $x_1^*, \dots, x_{n_p}^*$, the value of $E(x_j^*)$ as defined
  in~(\ref{errorj}).}
\label{fig3}
\end{figure}

Table~\ref{tab1b} shows the performance of Coordinate Descent, the
Spectral Projected Gradient (SPG)
method~\cite{bmr,bmr2,bmrreview,bmrenc}, and
Gencan~\cite{bmgencan,bmcomper}. In all cases, the initial point given
by the Fang-O'Leary technique was used. Since problem~(\ref{MDGP}) is
unconstrained, applying SPG corresponds to applying the Spectral
Gradient methods as proposed in~\cite{ghr}; while applying Gencan
corresponds to applying a line search Newton's method as considered
in~\cite{fang}. All three methods used as stopping criterion $f(x^k)
\leq \ftarget := 10^{-10}$. In addition, SPG and Gencan also stopped
if $\| \nabla f(x^k) \|_{\infty} \leq \varepsilon_{\mathrm{opt}} :=
10^{-8}$. For all three methods the table shows the number of
iterations (\#iter), the CPU time in seconds (Time), the value of the
objective function at the final iterate ($f(x^*)$), and the error with
respect to the known solution ($E(x^*)$). In addition, the table
shows, for the coordinate descent method the number of evaluations of
$\underline{f}$; while it shows for the other two methods, the number
of evaluations of $f$ and $\| \nabla f(x^*) \|_{\infty}$. In the
table, highlighted figures in column $f(x^*)$ are the ones that
correspond to local minimizers. Highlighted figures in column $E(x^*)$
correspond to final iterates that are far from the known solution. In
most cases, this fact is associated with having found a local
minimizer. However, in some cases, it corresponds to an alternative
global minimizer. We may observe that coordinate descent stands out as
the only method to have found a global minimizer in all the eighteen
considered instances. Figures~\ref{fig2} and~\ref{fig3} illustrate
three molecules in which the coordinate descent method found a global
solution while SPG and Gencan found local non-global minimizers. It is
worth mentioning that the numerical experiments reported
in~\cite{artacho} show that the Douglas-Rachford method, that requires
an SVD decomposition of a $n_p \times n_p$ matrix per iteration, with
a limit of $5{,}000$ iterations, was able to reconstruct the two
smallest molecules (1PTQ and 1HOE) only. As reported in
\cite{artacho}, the reconstruction of molecules 1LFB and 1PHT was
``satisfactory''; while the reconstruction of molecules 1POA and 1AX8
was ``poor''.

\begin{table}[ht!]
\begin{center}
\resizebox{\textwidth}{!}{
\begin{tabular}{|c|c|c|rrrcc|rrrccc|rrrccc|}
\hline
& & \multirow{2}{*}{Molecule} &  
\multicolumn{5}{c|}{Coordinate descent} &
\multicolumn{6}{c|}{Spectral Projected Gradient} &
\multicolumn{6}{c|}{Gencan} \\
\cline{4-20}
& & &
\multicolumn{1}{c}{\#iter} &  
\multicolumn{1}{c}{\#$\underline{f}$} &  
\multicolumn{1}{c}{Time} &  
\multicolumn{1}{c}{$f(x^*)$} &
\multicolumn{1}{c|}{$E(x^*)$} &
\multicolumn{1}{c}{\#iter} &  
\multicolumn{1}{c}{\#$f$} &  
\multicolumn{1}{c}{Time} &  
\multicolumn{1}{c}{$f(x^*)$} &
\multicolumn{1}{c}{$\| \nabla f(x^*) \|_{\infty}$} &
\multicolumn{1}{c|}{$E(x^*)$} &
\multicolumn{1}{c}{\#iter} &  
\multicolumn{1}{c}{\#$f$} &  
\multicolumn{1}{c}{Time} &  
\multicolumn{1}{c}{$f(x^*)$} &
\multicolumn{1}{c}{$\| \nabla f(x^*) \|_{\infty}$} &
\multicolumn{1}{c|}{$E(x^*)$} \\
\hline
\hline
\multirow{18}{*}{\rotatebox{90}{
\begin{tabular}{c}
Points may correspond to protein\\
atoms (ATOM) only or to protein atoms\\
plus atoms in small molecules (HETATM)
\end{tabular}}}
 & \multirow{12}{*}{\rotatebox{90}{ATOM only}}
 & 1ptq &     57,671 &     57,686 &  0.13 & 9.99e-11 & 1.66e-06 &    333 &    334 & 0.05 &   1.12e-11 & 3.32e-07 &   2.60e-06 &   9 &  13 &   0.42 &   5.82e-13 & 1.45e-07 &   2.58e-06 \\
&& 1hoe &    135,886 &    135,907 &  0.30 & 9.99e-11 & 2.36e-06 &    126 &    128 & 0.03 &   8.52e-11 & 4.39e-07 &   3.62e-06 &   7 &  10 &   0.62 &   5.49e-11 & 1.26e-06 &   7.25e-06 \\
&& 1lfb &    811,486 &    811,613 &  1.79 & 9.99e-11 & 7.56e-06 &    738 &    755 & 0.19 &   1.77e-11 & 5.19e-07 &   1.11e-06 &  13 &  20 &   1.15 &   5.96e-11 & 1.05e-06 &   4.27e-06 \\
&& 1pht &    786,655 &    786,831 &  1.99 & 9.99e-11 & 1.79e-05 &  5,945 &  6,856 & 2.50 &\h{2.85e-02}& 5.04e-09 &\h{2.08e-01}& 127 & 340 &  28.67 &\h{2.85e-02}& 9.10e-07 &\h{2.08e-01}\\
&& 1poa &    704,652 &    704,762 &  1.59 & 9.99e-11 & 9.79e-06 &  7,367 &  8,716 & 3.00 &   9.95e-11 & 3.00e-08 &   5.84e-04 &  18 &  19 &   2.71 &   1.79e-11 & 3.68e-07 &   2.32e-04 \\
&& 2igg &    484,388 &    484,473 &  1.56 & 9.99e-11 & 9.12e-06 &    304 &    305 & 0.21 &   8.79e-11 & 3.18e-07 &   3.43e-06 &  11 &  22 &   4.58 &   2.45e-13 & 1.24e-08 &   2.58e-07 \\
&& 1ax8 &    353,820 &    353,895 &  0.80 & 9.99e-11 & 2.54e-06 &    325 &    326 & 0.14 &   8.11e-11 & 1.74e-07 &   2.02e-05 &  14 &  20 &   3.59 &   7.98e-13 & 1.97e-08 &   2.14e-06 \\
&& 1rml &    340,528 &    340,586 &  1.24 & 9.99e-11 & 4.07e-06 &    236 &    238 & 0.39 &   1.21e-12 & 9.05e-08 &   1.46e-06 &   9 &  10 &  30.48 &   3.77e-13 & 6.47e-08 &   1.23e-06 \\
&& 1ak6 & 15,138,479 & 15,138,810 &229.85 & 9.99e-11 & 1.35e-05 &  1,662 &  1,755 & 4.06 &\h{5.18e-02}& 9.92e-09 &\h{1.97e-01}& 166 & 421 & 953.67 &\h{5.18e-02}& 7.47e-09 &\h{1.97e-01}\\
&& 1a24 &  2,840,577 &  2,840,834 &  9.87 & 9.99e-11 & 1.39e-05 &    322 &    325 & 0.74 &   8.10e-11 & 5.88e-08 &   1.15e-05 &  19 &  48 &  74.23 &   8.76e-12 & 1.76e-07 &   7.89e-06 \\
&& 3msp & 12,873,352 & 12,874,426 & 42.40 & 9.99e-11 & 1.61e-05 &    672 &    688 & 1.93 &   1.41e-10 & 8.47e-09 &   1.86e-05 &  31 &  55 & 138.70 &   1.24e-11 & 1.56e-08 &   5.46e-06 \\
&& 3eza & 17,122,466 & 17,123,479 & 58.89 & 9.99e-11 & 1.03e-05 &    580 &    586 & 2.26 &   4.43e-10 & 9.66e-09 &   2.11e-05 &  23 &  51 & 224.11 &   1.24e-10 & 3.08e-09 &   1.18e-05 \\
\cline{2-20}
 & \multirow{6}{*}{\rotatebox{90}{{\scriptsize ATOM+HETATM}}}
 & 1ptq &     57,640 &     57,659 &  0.13 & 9.99e-11 & 1.61e-06 &    334 &    335 & 0.05 &   8.45e-11 & 2.81e-07 &   3.60e-05 &  10 &  15 &   0.45 &   3.83e-17 & 9.27e-10 &   2.24e-08 \\
&& 1hoe &    129,571 &    129,590 &  0.31 & 9.99e-11 & 2.25e-06 &    143 &    144 & 0.04 &   3.48e-11 & 3.47e-07 &   2.82e-06 &   8 &  11 &   0.76 &   8.61e-18 & 4.80e-10 &   3.14e-09 \\
&& 1pht &    946,496 &    946,610 &  8.26 & 9.99e-11 & 1.60e-05 &  1,541 &  1,608 & 0.78 &\h{1.61e-05}& 9.85e-09 &\h{1.04e-01}&  21 &  31 &   8.48 &\h{1.61e-05}& 1.49e-09 &\h{1.04e-01}\\
&& 1poa &    409,610 &    409,655 &  0.93 & 9.99e-11 & 5.20e-05 & 12,996 & 15,710 & 6.43 &   1.84e-10 & 9.99e-09 &\h{1.57e-02}&  15 &  18 &   4.31 &   1.93e-11 & 4.49e-07 &\h{1.38e-02}\\
&& 1ax8 &    308,962 &    309,026 &  0.70 & 9.99e-11 & 2.18e-06 &    148 &    149 & 0.07 &   9.43e-11 & 2.25e-07 &   1.04e-05 &   8 &  11 &   3.18 &   1.54e-12 & 2.35e-07 &   6.49e-07 \\
&& 1rml &    344,977 &    345,021 &  1.28 & 9.99e-11 & 4.01e-06 &    305 &    307 & 0.52 &   9.48e-11 & 1.55e-07 &   4.36e-05 &  10 &  11 &  36.41 &   9.82e-17 & 5.53e-10 &   4.35e-08 \\
\hline
\end{tabular}}
\end{center}
\caption{Performance of Coordinate Descent, SPG, and Gencan applied to
  the instances of problem~(\ref{MDGP}) built with the molecules
  considered in~\cite{artacho} or~\cite{fang}.}
\label{tab1b}
\end{table}

At this point the following question arises: how does solving the
subproblems with cubically-regularized second-order models affect the
performance of the CD method? To answer this question, we solved the
same 18 problems tackling the subproblems with
quadratically-regularized linear models. This means that, to
approximately solve~(\ref{ell}) with Algorithm~\ref{preliminar}.1, we
considered $p=1$. In other words, instead
of~(\ref{modelellmenosuno},\ref{modelell}), \textbf{(a)} we considered
the first-order Taylor expansion of~$\underline{f}$ at $\bar x =
x_{\ell(k)}^k \in \R^d$ given by $M_{\bar x}(z) := \underline{f}(\bar
x) + \nabla \underline{f}(\bar x)^T (z - \bar x)$, and \textbf{(b)} we
computed $\xtrial$ as the global minimizer of
\begin{equation} \label{queseyo}
  \Minimize_{z \in \R^d} M_{\bar x}(z) + \sigma \| z - \bar x \|^2.
\end{equation}
Since (\ref{queseyo}) has no solution when $\nabla \underline{f}(\bar
x) \neq 0$ and $\sigma=0$, we skip the case $\sigma=0$ by substituting
$\sigma \leftarrow 0$ with $\sigma \leftarrow \sigma_{\min}$ at Step~1
of Algorithm~\ref{preliminar}.1. Apart from this, the settings for the
case $p=1$ were identical to those already described for the case
$p=2$. Table~\ref{tab1c} shows the results. The numbers in the table
show that the method found a global solution in all instances, a
feature shared with its counterpart with $p=2$. (Only in one instance
an alternative global minimizer was found.) The numbers in the table
also show that, on average, the method does 1.0001 function
evaluations per iteration when $p=2$, while that same amount is 1.5000
when $p=1$. This means that, on the one hand, in the case $p=1$, half
of the times the solution of the regularized model is discarded for
not satisfying the descent condition and the regularization parameter
must be increased. On the other hand, this situation is extremely rare
(once every ten thousand iterations) when $p=2$.  Moreover, the method
with $p=1$ uses, on average, 26 times more iterations, 39 times more
function evaluations and 22 times more time than the case $p=2$. The
conclusion is that using quadratic models with cubic regularization
whose global solution can be calculated using the method introduced
in~\cite{bmbunch}, greatly improves the performance of the proposed
method.

\begin{table}[ht!]
\begin{center}
\begin{tabular}{|c|c|rrrcc|}
\hline
& \multirow{2}{*}{Molecule} & \multicolumn{5}{c|}{Coordinate descent with $p=1$}\\
\cline{3-7}
& &
\multicolumn{1}{c}{\#iter} &  
\multicolumn{1}{c}{\#$\underline{f}$} &  
\multicolumn{1}{c}{Time} &  
\multicolumn{1}{c}{$f(x^*)$} &
\multicolumn{1}{c|}{$E(x^*)$}\\
\hline
\hline
\multirow{12}{*}{\rotatebox{90}{ATOM only}}
& 1ptq &   1,115,103 &   1,672,658 &   1.97 & 9.99e-11 & 3.98e-05\\
& 1hoe &   1,862,919 &   2,794,382 &   3.10 & 9.99e-11 & 2.04e-06\\
& 1lfb &   7,059,295 &  10,588,949 &   9.50 & 9.99e-11 & 7.52e-06\\
& 1pht &  21,056,147 &  31,584,233 &  38.52 & 9.99e-11 & 2.73e-04\\
& 1poa & 107,284,920 & 160,927,383 & 105.32 & 9.99e-11 & 5.85e-04\\
& 2igg &   3,602,906 &   5,404,362 &  30.38 & 9.99e-11 & 8.43e-06\\
& 1ax8 &   2,059,932 &   3,089,901 &   4.56 & 9.99e-11 & 4.05e-06\\
& 1rml &   9,186,681 &  13,780,025 &  98.21 & 9.99e-11 & 3.23e-06\\
& 1ak6 & 160,256,027 & 240,384,044 & 465.43 & 9.99e-11 & 1.62e-05\\
& 1a24 &  79,008,635 & 118,512,956 & 238.59 & 9.99e-11 & 1.40e-05\\
& 3msp & 236,434,288 & 354,651,435 & 458.60 & 9.99e-11 & 1.61e-05\\
& 3eza &  34,134,889 &  51,202,336 & 216.98 & 9.99e-11 & 1.07e-05\\
\hline
\multirow{6}{*}{\rotatebox{90}{{\scriptsize ATOM+HETATM}}}
& 1ptq &   1,466,069 &    2199,107 &   2.30 & 9.97e-11 & 4.26e-05\\
& 1hoe &     683,645 &    1025,474 &   3.53 & 9.99e-11 & 1.98e-06\\
& 1pht &  16,112,116 &   24168,177 &  23.29 & 9.99e-11 & 3.38e-05\\
& 1poa &  32,496,533 &   48744,802 &  34.40 & 9.99e-11 & \h{1.46e-02}\\
& 1ax8 &   6,829,569 &   10244,357 &   9.37 & 9.99e-11 & 2.47e-06\\
& 1rml &   2,302,956 &    3454,437 &  87.64 & 9.99e-11 & 1.99e-06\\
\hline
\end{tabular}
\end{center}
\caption{Performance of Coordinate Descent with $p=1$,
  i.e.\ considering quadratically-regularized linear models for
  solving subproblems, applied to the same instances already shown in
  Table~\ref{tab1b}.}
\label{tab1c}
\end{table}

Another natural question that arises is whether the tendency of the
coordinate descent method in finding global minimizers could be
observed in a larger set of instances. To check this hypothesis, we
downloaded 64 additional \textit{random} molecules with no more than
$6{,}000$ atoms from the ones that were uploaded in 2020; 56 of which
have, other than protein atoms, atoms in small molecules. However
there were 19 molecules for which, considering protein atoms only or
protein atoms plus atoms in small molecules, the graph associated with
the incomplete Euclidean matrix obtained by eliminating distances
larger than 6 Angstroms is disconnected. Therefore, we were left with
45 and 37 molecules in each set, totalizing 82 new
instances. Table~\ref{tab2} shows the performance of Coordinate
Descent and SPG when applied to the 45 instances that consider protein
atoms only; while Table~\ref{tab3} shows the performance of both
methods when applied to the 37 instances that consider protein atoms
plus atoms in small molecules. In the 45 instances in
Table~\ref{tab2}, Coordinate Descent found 37 global minimizers; while
SPG found 30 global minimizers.
%This means that Coordinate Descent found 23\% more global minimizers than SPG.
In the 37 instances in Table~\ref{tab3}, Coordinate Descent found 30
global minimizers; while SPG found 26 global minimizers.
%This means that Coordinate Descent found 15\% more global minimizers than SPG.

\begin{table}[ht!]
\begin{center}
\resizebox{\textwidth}{!}{
\begin{tabular}{|c|rrrr|rrrcc|rrrccc|}
\hline
\multirow{2}{*}{Molecule} &  
\multicolumn{1}{c}{\multirow{2}{*}{$n$}} &  
\multicolumn{1}{c}{\multirow{2}{*}{$n_p$}} &  
\multicolumn{1}{c}{\multirow{2}{*}{$|S|$}} &  
\multicolumn{1}{c|}{\multirow{2}{*}{Time $x^0$}} &  
\multicolumn{5}{c|}{Coordinate descent} &
\multicolumn{6}{c|}{Spectral Projected Gradient} \\
\cline{6-16}
 & & & & &
\multicolumn{1}{c}{\#iter} &  
\multicolumn{1}{c}{\#$\underline{f}$} &  
\multicolumn{1}{c}{Time} &  
\multicolumn{1}{c}{$f(x^*)$} &
\multicolumn{1}{c|}{$E(x^*)$} &
\multicolumn{1}{c}{\#iter} &  
\multicolumn{1}{c}{\#$f$} &  
\multicolumn{1}{c}{Time} &  
\multicolumn{1}{c}{$f(x^*)$} &
\multicolumn{1}{c}{$\| \nabla f(x^*) \|_{\infty}$} &
\multicolumn{1}{c|}{$E(x^*)$} \\
\hline
\hline
6kbq &  7,554 & 2,518 &  98,650  1.56 &  40.49 &   10,543,341 &   10,543,995 &   36.22 &   9.99e-11 &   2.65e-05 &    810 &    911 &  0.96 &\h{2.65e-02}& 7.65e-09 &\h{1.87e-01}\\%ONLY CD
6kc2 &  7,530 & 2,510 &  98,480  1.56 &  39.97 &    7,682,589 &    7,683,021 &   20.93 &   9.99e-11 &   2.62e-05 &    954 &    984 &  1.11 &   4.48e-12 & 3.43e-08 &   1.29e-05 \\%B
6khu &  3,147 & 1,049 &  39,372  3.58 &   3.09 &     ,316,876 &      316,959 &    0.83 &   9.99e-11 &   6.12e-06 &    411 &    413 &  0.19 &   5.66e-11 & 1.05e-06 &   7.83e-06 \\%B
6kir &  3,861 & 1,287 &  48,596  2.94 &   5.58 &    7,574,444 &    7,574,767 &   54.37 &\h{1.60e-03}&\h{1.50e-01}&    871 &    895 &  0.49 &   5.27e-11 & 1.01e-07 &   3.82e-05 \\%ONLY SPG
6kk9 & 17,538 & 5,846 & 218,662  0.64 & 479.82 &  215,549,505 &  215,561,546 &  712.59 &\h{9.29e-01}&\h{1.19e+00}&  6,262 &  7,030 & 17.33 &\h{9.09e-01}& 9.83e-09 &\h{1.15e+00}\\%NONE
6kki &  8,193 & 2,731 & 109,472  1.47 &  53.06 &    2,113,286 &    2,113,472 &    5.68 &   9.99e-11 &   5.04e-06 &    432 &    434 &  0.56 &   9.34e-11 & 1.33e-07 &   2.99e-05 \\%B
6kkj &  7,866 & 2,622 & 105,642  1.54 &  47.70 &    4,068,836 &    4,070,082 &   10.88 &   9.99e-11 &   5.18e-06 &  5,868 &  6,787 &  7.70 &   3.24e-10 & 9.99e-09 &   1.58e-03 \\%B
6kkl &  8,262 & 2,754 & 112,646  1.49 &  54.76 &    6,064263  &    6,064,648 &   16.61 &   9.99e-11 &   5.68e-06 &    592 &    604 &  0.78 &   9.47e-11 & 5.99e-08 &   6.48e-06 \\%B
6kkv & 14,607 & 4,869 & 357,016  1.51 & 290.74 &   31,543,917 &   31,545,196 &  128.60 &   9.99e-11 &   1.86e-05 &  1,095 &  1,140 &  4.65 &   2.88e-11 & 6.93e-08 &   2.00e-06 \\%B
6kx0 &  9,885 & 3,295 & 128,930  1.19 &  91.69 &   38,237,232 &   38,239,417 &  196.30 &\h{1.57e-01}&\h{1.02e+00}&  1,428 &  1,546 &  2.26 &\h{1.57e-01}& 9.87e-09 &\h{1.02e+00}\\%NONE
6kys &  5,442 & 1,814 &  70,296  2.14 &  16.13 &      984,382 &      984,559 &    2.49 &   9.99e-11 &   4.00e-06 &    359 &    360 &  0.29 &   8.44e-11 & 1.22e-07 &   2.95e-05 \\%B
6l29 &  5,427 & 1,809 &  70,480  2.15 &  16.68 &      965,426 &      965,620 &    2.49 &   9.99e-11 &   4.51e-06 &    433 &    436 &  0.36 &   8.43e-11 & 7.93e-07 &   8.26e-06 \\%B
6l2a &  2,463 &   821 &  29,876  4.44 &   1.55 &      411,298 &      411,381 &    1.04 &   9.99e-11 &   5.55e-06 &  1,612 &  1,699 &  0.59 &   9.95e-11 & 7.47e-08 &   6.46e-05 \\%B
6laf &  8,442 & 2,814 & 105,164  1.33 &  58.06 &   31,413,895 &   31,418,665 &  104.66 &   9.99e-11 &   2.51e-05 &  7,706 &  9,068 & 10.22 &\h{1.53e-02}& 9.76e-09 &\h{2.74e-01}\\%ONLY CD
6li7 &  7,515 & 2,505 &  98,374  1.57 &  40.21 &   10,533,446 &   10,534,072 &   37.27 &   9.99e-11 &   2.64e-05 &    698 &    713 &  0.83 &   8.73e-11 & 4.32e-07 &   1.54e-05 \\%B
6lik &  7,515 & 2,505 &  97,354  1.55 &  39.67 &    7,579,646 &    7,580,087 &   19.88 &   9.99e-11 &   2.63e-05 &    727 &    749 &  0.85 &   8.89e-11 & 2.61e-07 &   6.42e-05 \\%B
6lty &  5,964 & 1,988 &  73,296  1.86 &  19.95 &    3,800,264 &    3,800,792 &    9.82 &   9.99e-11 &   7.86e-06 &    672 &    683 &  0.58 &   9.76e-11 & 9.54e-08 &   1.78e-05 \\%B
6ltz &  3,432 & 1,144 &  41,230  3.15 &   3.97 &     ,975,177 &      975,323 &    2.48 &   9.99e-11 &   3.13e-06 &  8,025 &  9,404 &  4.08 &   1.27e-10 & 9.94e-09 &   9.66e-04 \\%B
6m37 &  9,999 & 3,333 & 125,394  1.13 &  91.51 &  130,035,728 &  130,040,179 &  422.49 &\h{5.07e-02}&\h{4.33e-01}&  4,936 &  5,479 &  7.63 &\h{5.06e-02}& 9.99e-09 &\h{4.30e-01}\\%NONE
6m5n &  5,781 & 1,927 &  73,890  1.99 &  18.48 &    8,183,149 &    8,183,657 &   27.54 &   9.99e-11 &   6.22e-06 &  6,419 &  7,461 &  5.81 &\h{3.06e-03}& 9.98e-09 &\h{2.74e-01}\\%ONLY CD
6m6j &    759 &   253 &  13,914 21.82 &   0.07 &       13,076 &       13,076 &    0.04 &   9.99e-11 &   1.84e-06 &     59 &     61 &  0.01 &   3.99e-12 & 1.99e-07 &   6.82e-07 \\%B
6m6k &    747 &   249 &  13,672 22.14 &   0.07 &       11,376 &       11,376 &    0.04 &   9.97e-11 &   1.29e-06 &     53 &     55 &  0.01 &   2.83e-11 & 5.59e-07 &   1.41e-06 \\%B
6pq0 &  6,087 & 2 029 &  84,104  2.04 &  22.54 &    2,346,681 &    2,346,842 &   31.31 &   9.99e-11 &   3.83e-06 &    440 &    442 &  0.43 &\h{2.39e-04}& 8.43e-09 &\h{1.83e-01}\\%ONLY CD
6pup &  4,035 & 1,345 &  50,372  2.79 &   6.57 &      544,519 &      544,617 &    1.42 &   9.99e-11 &   4.87e-06 &    409 &    410 &  0.24 &   2.52e-11 & 7.99e-08 &   2.38e-05 \\%B
6pxf &  4,944 & 1,648 &  64,944  2.39 &  11.70 &      811,879 &      812,002 &    2.13 &   9.99e-11 &   2.77e-05 &  4,907 &  5,583 &  3.88 &   1.81e-10 & 9.41e-09 &   2.91e-03 \\%B
6q08 &    741 &   247 &  13,232 21.78 &   0.07 &       22,723 &       22,744 &    0.07 &   9.99e-11 &   7.19e-06 &    173 &    176 &  0.03 &   5.75e-11 & 1.85e-06 &   3.82e-06 \\%B
6sx6 &  2,340 &   780 &  44,862  7.38 &   1.28 &      143,230 &      143,260 &    0.49 &   9.99e-11 &   2.50e-06 &    136 &    137 &  0.07 &   5.51e-12 & 1.80e-07 &   2.60e-07 \\%B
6syk &  2,718 &   906 &  54,052  6.59 &   1.99 &      156,648 &      156,707 &    0.55 &   9.99e-11 &   2.50e-06 &    585 &    592 &  0.36 &   9.41e-12 & 2.34e-07 &   1.42e-05 \\%B
6t1z &  8,943 & 2,981 & 119,580  1.35 &  65.75 &   14,818,888 &   14,819,492 &   40.30 &   9.99e-11 &   8.88e-06 &  3,247 &  3,539 &  4.71 &   9.94e-11 & 1.37e-08 &   2.70e-04 \\%B
6tad &  4,362 & 1,454 &  58,736  2.78 &   8.07 &    2,385,243 &    2,385,423 &    6.35 &   9.99e-11 &   6.84e-06 &    632 &    643 &  0.44 &   7.87e-11 & 3.81e-07 &   5.70e-06 \\%B
6twe &  7,902 & 2,634 & 163,598  2.36 &   0.18 &    5,773,930 &    5,774,426 &   20.85 &   9.99e-11 &   8.20e-06 &    598 &    612 &  1.16 &   7.04e-11 & 1.66e-07 &   6.93e-06 \\%B
6ubh &  9,009 & 3,003 & 113,766  1.26 &  69.74 &    8,782,138 &    8,782,736 &   22.99 &   9.99e-11 &   1.64e-05 &  3,565 &  4,019 &  4.95 &   1.23e-10 & 9.86e-09 &   1.63e-03 \\%B
6ucd &  8,199 & 2,733 &  97,910  1.31 &  52.10 &   86,618,396 &   86,621,926 &  278.99 &\h{1.35e+00}&\h{1.88e+00}&  9,391 & 11,180 & 11.51 &\h{9.92e-01}& 9.29e-09 &\h{1.10e+00}\\%NONE
6veh &  7,431 & 2,477 & 182,676  2.98 &  42.54 &    4,087,012 &    4,087,391 &   16.46 &   9.99e-11 &   1.43e-05 &    658 &    662 &  1.37 &   3.50e-11 & 1.38e-07 &   5.83e-06 \\%B
6vk2 &  4,704 & 1,568 & 108,520  4.42 &  11.94 &      463,141 &      463,360 &    1.79 &   9.99e-11 &   9.46e-06 &  3,069 &  3,419 &  3.95 &   1.40e-10 & 9.62e-09 &   2.22e-03 \\%B
6vnz &  1,392 &   464 &  25,364 11.81 &   0.32 &       88,484 &       88,530 &    0.29 &   9.99e-11 &   3.40e-06 &    206 &    208 &  0.06 &   6.02e-11 & 2.16e-07 &   4.36e-06 \\%B
6vv6 &  7,464 & 2,488 &  92,878  1.50 &  45.42 &    7,583,602 &    7,584,292 &   33.37 &   9.99e-11 &   4.04e-06 &  1,846 &  1,943 &  2.05 &\h{6.23e-02}& 1.00e-08 &\h{2.23e-01}\\%ONLY CD
6vv7 &  7,452 & 2,484 &  93,102  1.51 &  44.46 &    7,692,189 &    7,692,886 &   41.06 &   9.99e-11 &   3.80e-06 &  1,981 &  2,233 &  2.26 &\h{1.07e-01}& 7.39e-09 &\h{2.23e-01}\\%ONLY CD
6vv9 &  7,452 & 2,484 &  92,506  1.50 &  44.29 &    6,251,096 &    6,251,759 &   29.79 &   9.99e-11 &   3.11e-06 &  1,022 &  1,040 &  1.14 &\h{7.95e-02}& 1.25e-09 &\h{2.17e-01}\\%ONLY CD
6wcr & 11,766 & 3,922 & 151,164  0.98 & 158.07 &  118,185,068 &  118,209,777 &  421.33 &\h{1.67e-01}&\h{3.74e-01}&  7,460 &  8,498 & 13.92 &\h{1.99e-01}& 9.96e-09 &\h{2.38e-01}\\%NONE
6yuc &  8,637 & 2,879 & 107,258  1.29 &  60.03 &   17,737,727 &   17,738,988 &  117.20 &\h{3.77e-01}&\h{6.86e-01}&  1,450 &  1,668 &  1.91 &\h{1.09e-04}& 9.93e-09 &\h{1.19e-01}\\%NONE
6z4c &  5,838 & 1,946 &  72,838  1.92 &  18.87 &    1,951,179 &    1,951,337 &    5.00 &   9.99e-11 &   2.72e-06 &    369 &    370 &  0.31 &   6.24e-11 & 2.43e-07 &   2.91e-06 \\%B
6zcm &  7,899 & 2,633 & 102,030  1.47 &  46.05 &   12,338,443 &   12,339,335 &   40.78 &   9.99e-11 &   2.53e-05 &  7,445 &  8,812 &  9.41 &   1.38e-10 & 9.93e-09 &   9.87e-04 \\%B
7ckj &  4,731 & 1,577 &  59,608  2.40 &  10.24 &    2,681,962 &    2,682,163 &   11.61 &   9.99e-11 &   3.27e-06 &    768 &    784 &  0.54 &\h{6.51e-07}& 5.85e-09 &\h{2.03e-01}\\%ONLY CD
7jjl &  9,690 & 3,230 & 125,976  1.21 &  83.34 &   77,851,519 &   77,854,837 &  303.65 &\h{4.82e-01}&\h{4.52e-01}&  4,830 &  5,451 &  7.49 &\h{7.94e-01}& 8.12e-09 &\h{4.73e-01}\\%NONE
\hline
\end{tabular}}
\end{center}
\caption{Performance of Coordinate Descent and SPG methods in the 46
  instances that consider protein atoms only.}
\label{tab2}
\end{table}

\begin{table}[ht!]
\begin{center}
\resizebox{\textwidth}{!}{
\begin{tabular}{|c|rrrr|rrrcc|rrrccc|}
\hline
\multirow{2}{*}{Molecule} &  
\multicolumn{1}{c}{\multirow{2}{*}{$n$}} &  
\multicolumn{1}{c}{\multirow{2}{*}{$n_p$}} &  
\multicolumn{1}{c}{\multirow{2}{*}{$|S|$}} &  
\multicolumn{1}{c|}{\multirow{2}{*}{Time $x^0$}} &  
\multicolumn{5}{c|}{Coordinate descent} &
\multicolumn{6}{c|}{Spectral Projected Gradient} \\
\cline{6-16}
 & & & & &
\multicolumn{1}{c}{\#iter} &  
\multicolumn{1}{c}{\#$\underline{f}$} &  
\multicolumn{1}{c}{Time} &  
\multicolumn{1}{c}{$f(x^*)$} &
\multicolumn{1}{c|}{$E(x^*)$} &
\multicolumn{1}{c}{\#iter} &  
\multicolumn{1}{c}{\#$f$} &  
\multicolumn{1}{c}{Time} &  
\multicolumn{1}{c}{$f(x^*)$} &
\multicolumn{1}{c}{$\| \nabla f(x^*) \|_{\infty}$} &
\multicolumn{1}{c|}{$E(x^*)$} \\
\hline
\hline
6kbq &  8,106 & 2,702 & 108,066  ( 1.48\%) &  50.64 &   6,476,123 &   6,476,435 &  17.57 &   9.99D-11 &   2.35D-05 &    571 &    580 &  0.71 &   9.97D-11 & 6.32D-08 &   1.85D-05 \\%B
6kc2 &  8,397 & 2,799 & 111,198  ( 1.42\%) &  56.89 &  15,115,487 &  15,118,534 &  50.50 &   9.99D-11 &\h{2.15D-01}&  1,102 &  1,121 &  1.42 &   4.03D-10 & 9.36D-09 &\h{2.19D-01}\\%B
6khu &  3,465 & 1,155 &  44,672  ( 3.35\%) &   4.24 &     334,415 &     334,473 &   0.87 &   9.99D-11 &   6.35D-06 &    490 &    500 &  0.25 &   6.20D-11 & 1.63D-06 &   7.45D-06 \\%B
6kir &  4,203 & 1,401 &  54,634  ( 2.79\%) &   7.42 &     560,278 &     560,406 &   1.49 &   9.99D-11 &   4.20D-06 &    450 &    454 &  0.28 &   9.28D-11 & 4.63D-07 &   6.28D-05 \\%B
6kk9 & 17,895 & 5,965 & 225,776  ( 0.63\%) & 515.33 & 202,467,161 & 202,477,953 & 691.77 &\h{4.00D-01}&\h{3.09D-01}&  5,084 &  5,637 & 13.87 &\h{3.82D-01}& 9.95D-09 &\h{3.70D-01}\\%NONE
6kki &  8,364 & 2,788 & 111,596  ( 1.44\%) &  56.23 &   2,049,210 &   2,049,408 &   5.51 &   9.99D-11 &   3.93D-06 &    592 &    594 &  0.76 &   9.92D-11 & 1.13D-07 &   5.68D-05 \\%B
6kkj &  7,992 & 2,664 & 106,472  ( 1.50\%) &  49.62 &  17,809,921 &  17,816,613 & 132.85 &\h{9.54D-04}&\h{2.81D-01}&  2,921 &  3,188 &  3.67 &\h{9.54D-04}& 9.12D-09 &\h{2.81D-01}\\%NONE
6kkl &  8,421 & 2,807 & 114,184  ( 1.45\%) &  58.35 &   6,082,607 &   6,082,924 &  16.57 &   9.99D-11 &   5.67D-06 &    521 &    537 &  0.70 &   6.88D-11 & 2.32D-07 &   6.82D-06 \\%B
6kkv & 14,820 & 4,940 & 364,756  ( 1.49\%) & 304.81 &  27,357,834 &  27,358,998 & 113.43 &   9.99D-11 &   1.80D-05 &  1,193 &  1,220 &  5.01 &   7.35D-10 & 7.99D-09 &   4.74D-05 \\%B
6kx0 & 10,044 & 3,348 & 131,624  ( 1.17\%) &  94.64 &  45,983,700 &  45,985,561 & 226.68 &\h{7.16D-02}&\h{2.47D-01}&  1,044 &  1,096 &  1.61 &\h{1.78D-01}& 9.27D-09 &\h{1.01D+00}\\%NONE
6kys &  5,886 & 1,962 &  77,926  ( 2.03\%) &  20.28 &     816,852 &     816,930 &   2.19 &   9.99D-11 &   3.60D-06 &  2,389 &  2,688 &  2.22 &   9.98D-11 & 3.11D-08 &   1.72D-04 \\%B
6l29 &  5,841 & 1,947 &  77,602  ( 2.05\%) &  19.25 &   1,807,474 &   1,807,572 &  18.75 &   9.99D-11 &\h{1.26D-01}&    589 &    592 &  0.52 &\h{5.31D-04}& 4.39D-09 &\h{1.60D-01}\\%ONLY CD
6l2a &  2,643 &   881 &  32,644  ( 4.21\%) &   1.82 &     787,753 &     787,881 &   4.57 &   9.99D-11 &   5.77D-06 &  1,769 &  1,892 &  0.67 &   1.79D-10 & 9.03D-09 &\h{1.24D-01}\\%B
6laf &  8,535 & 2,845 & 106,084  ( 1.31\%) &  57.40 &  38,880,745 &  38,890,853 & 174.71 &\h{1.57D-02}&\h{3.40D-01}& 10,824 & 12,773 & 13.95 &\h{1.71D-02}& 9.71D-09 &\h{3.41D-01}\\%NONE
6li7 &  8,520 & 2,840 & 114,120  ( 1.42\%) &  57.40 &   6,141,620 &   6,141,906 &  16.33 &   9.99D-11 &   2.29D-05 &    443 &    453 &  0.58 &   9.37D-11 & 5.62D-08 &   1.50D-05 \\%B
6lik &  8,208 & 2,736 & 108,364  ( 1.45\%) &  51.66 &  12,686,358 &  12,686,897 &  52.81 &   9.99D-11 &   2.48D-05 &  7,235 &  8,363 &  9.40 &\h{2.99D-03}& 9.91D-09 &\h{1.81D-01}\\%ONLY CD
6ltz &  3,798 & 1,266 &  47,246  ( 2.95\%) &   5.62 &     767,531 &     767,612 &   1.93 &   9.99D-11 &   7.56D-06 &    985 &  1,002 &  0.54 &   9.73D-11 & 5.68D-08 &   1.28D-04 \\%B
6m5n &  6,687 & 2,229 &  88,996  ( 1.79\%) &  29.41 &   2,380,904 &   2,381,175 &   6.33 &   9.99D-11 &\h{1.41D-01}& 14,921 & 17,730 & 15.90 &   6.98D-10 & 9.55D-09 &\h{1.42D-01}\\%B
6m6j &    840 &   280 &  15,602  (19.97\%) &   0.09 &      18,319 &      18,319 &   0.06 &   9.99D-11 &   1.86D-06 &     76 &     78 &  0.01 &   9.94D-11 & 3.25D-07 &   2.06D-06 \\%B
6m6k &    828 &   276 &  15,388  (20.27\%) &   0.09 &      15,964 &      15,964 &   0.05 &   9.99D-11 &   1.63D-06 &     74 &     76 &  0.01 &   1.00D-11 & 3.36D-07 &   7.45D-07 \\%B
6pq0 &  6,360 & 2,120 &  89,018  ( 1.98\%) &   0.12 &   1,086,898 &   1,086,994 &  12.19 &   9.99D-11 &   3.98D-06 &  3,540 &  4,285 &  3.81 &   9.00D-11 & 1.21D-07 &   6.60D-05 \\%B
6pup &  4,272 & 1,424 &  55,442  ( 2.74\%) &   7.63 &     593,071 &     593,135 &   1.55 &   9.99D-11 &   3.60D-06 &    363 &    372 &  0.23 &   3.28D-11 & 6.55D-07 &   2.12D-06 \\%B
6pxf &  5,661 & 1,887 &  76,154  ( 2.14\%) &  17.47 &   6,753,157 &   6,775,873 &  25.50 &   9.97D-11 &\h{1.10D-01}&  4,541 &  5,069 &  4.07 &\h{1.04D-05}& 9.49D-09 &\h{1.16D-01}\\%ONLY CD
6t1z &  9,492 & 3,164 & 127,488  ( 1.27\%) &  79.93 &   3,427,330 &   3,427,847 &   9.17 &   9.99D-11 &\h{1.14D-02}&  7,688 &  9,017 & 11.70 &   2.94D-10 & 9.82D-09 &\h{1.23D-02}\\%B
6tad &  5,172 & 1,724 &  71,178  ( 2.40\%) &  13.52 &   1,050,635 &   1,050,702 &   2.89 &   9.99D-11 &   3.64D-06 &    273 &    275 &  0.22 &   9.83D-11 & 3.61D-08 &   3.95D-06 \\%B
6twe &  7,905 & 2,635 & 163,694  ( 2.36\%) &  45.84 &   5,768,235 &   5,768,740 &  20.83 &   9.99D-11 &   8.19D-06 &    637 &    643 &  1.17 &   9.14D-11 & 1.78D-07 &   8.01D-06 \\%B
6ubh &  9,999 & 3,333 & 130,348  ( 1.17\%) &  91.62 &  22,598,823 &  22,602,027 & 107.35 &   9.99D-11 &\h{3.60D-01}&  1,234 &  1,279 &  1.90 &   6.89D-11 & 3.25D-07 &\h{3.63D-01}\\%B
6veh &  7,566 & 2,522 & 186,586  ( 2.93\%) &  42.40 &   4,241,879 &   4,242,258 &  17.14 &   9.99D-11 &   1.52D-05 &    381 &    384 &  0.80 &   8.74D-11 & 2.80D-07 &   4.54D-05 \\%B
6vv6 &  8,169 & 2,723 & 107,072  ( 1.44\%) &  52.54 &   2,824,276 &   2,824,520 &   7.62 &   9.99D-11 &   9.75D-06 &    391 &    395 &  0.49 &   7.43D-11 & 8.50D-08 &   3.05D-06 \\%B
6vv7 &  8,247 & 2,749 & 109,326  ( 1.45\%) &  52.54 &   3,807,178 &   3,807,395 &  10.23 &   9.99D-11 &   2.33D-06 &  4,434 &  5,002 &  5.83 &   5.68D-11 & 3.59D-07 &   2.89D-04 \\%B
6vv9 &  8,250 & 2,750 & 108,784  ( 1.44\%) &  52.41 &   3,494,562 &   3,494,812 &   9.30 &   9.99D-11 &   8.41D-06 &    404 &    409 &  0.51 &   6.69D-11 & 2.00D-07 &   7.49D-06 \\%B
6wcr & 12,225 & 4,075 & 157,032  ( 0.95\%) & 167.40 & 229,852,446 & 229,885,273 & 721.26 &\h{9.77D-05}&\h{8.68D-02}&  6,972 &  7,906 & 13.17 &\h{8.93D-02}& 9.97D-09 &\h{1.97D-01}\\%NONE
6yuc &  8,640 & 2,880 & 107,366  ( 1.29\%) &  59.34 &  50,987,720 &  50,988,983 & 203.00 &\h{1.09D-04}&\h{1.19D-01}&  1,879 &  1,942 &  2.37 &\h{1.09D-04}& 7.65D-09 &\h{1.19D-01}\\%NONE
6z4c &  5,994 & 1,998 &  76,300  ( 1.91\%) &  20.35 &   1,697,728 &   1,697,882 &   4.45 &   9.99D-11 &   2.95D-06 &    270 &    271 &  0.24 &   9.71D-11 & 2.20D-07 &   3.12D-06 \\%B
6zcm &  9,696 & 3,232 & 135,448  ( 1.30\%) &  83.59 &   6,978,627 &   6,978,943 &  19.30 &   9.99D-11 &   7.42D-06 &  2,131 &  2,288 &  3.43 &   9.67D-11 & 5.81D-08 &   1.55D-04 \\%B
7ckj &  5,199 & 1,733 &  68,624  ( 2.29\%) &  13.55 &   1,630,610 &   1,630,694 &  10.52 &   9.99D-11 &   4.11D-06 &    541 &    549 &  0.43 &\h{2.14D-05}& 9.60D-09 &\h{1.96D-01}\\%ONLY CD
7jjl &  9,693 & 3,231 & 126,082  ( 1.21\%) &  83.84 & 102,233,476 & 102,237,339 & 369.96 &\h{2.46D-01}&\h{3.25D-01}&  5,332 &  6,097 &  8.12 &\h{7.94D-01}& 8.91D-09 &\h{4.73D-01}\\%NONE
\hline                        
\end{tabular}}                
\end{center}
\caption{Performance of Coordinate Descent and SPG methods in the 37
  instances that consider protein atoms plus atoms in small
  molecules.}
\label{tab3}
\end{table}

\section{Conclusions} \label{concl}

Methods based on high-order models for optimization are difficult to
implement due to the necessity of computing and storing high-order
derivatives and the complexity of solving the subproblems. These
difficulties are not so serious if the subproblems are
low-dimensional, which is the most frequent situation in the case of
CD methods. In the extreme case, in which one solves only univariate
problems, the number of high-order partial derivatives that are
necessary is a small multiple of the number of variables. Therefore,
the theory that shows that CD algorithms with high-order models enjoy
good convergence and complexity properties seems to be useful to
support the efficiency of practical implementations. In this context,
higher-order techniques allow to escape from attraction points that
tend to satisfy lower-order optimality conditions; see \cite{mr2014}.

Sometimes the fulfillment of a necessary high-order optimality
condition can be expressed as fulfillment of $\Phi(x) = 0$, where
$\Phi$ is a continuous nonnegative function. In this case, it makes
sense to say that $\Phi(x) \leq \varepsilon$ is an approximate
high-order optimality condition. Moreover, instead of requiring
globality for the solution to the regularized model-based
subproblem~(\ref{subproblem}), we may require only that $\Phi(x^{k+1})
\to 0$ when $k \to +\infty$, where $\Phi$ corresponds to the
high-order optimality condition of~(\ref{subproblem}). Careful choices
of $\Phi$ and the subproblems' stopping criterion may give rise to
complexity results associated with the attainment of these high-order
optimality conditions, see~\cite{cagt5,cagt6,cagt7}. This will be the
subject of future research.

In this paper the defined algorithms were applied to the
identification of proteins under NMR data. Moreover, we extended the
CD approach to the computation of a suitable initial approximation
that avoids, in many cases, the convergence to local non-global
minimizers. Our choice of the most adequate parameter $p$, that
defines the approximating models, and the strategy for choosing the
groups of variables were dictated by theoretical considerations
discussed in Section~\ref{discussion} and by the specific
characteristics of the problem. Our computing results are fully
reproducible and the codes are available in
\url{http://www.ime.usp.br/~egbirgin/}.

In future works we will apply the new CD techniques to the case in
which data uncertainty is present and outliers are likely to
occur. Possible improvements also include the choice of different
models at each iteration or at each group of variables with the aim of
making a better use of current information.\\

\noindent
\textbf{Data availability:} The datasets generated during and/or
analyzed during the current study are available in the corresponding
author web page, \url{http://www.ime.usp.br/~egbirgin/}.

\end{document}